\documentclass[conference, 10pt]{IEEEtran}
\IEEEoverridecommandlockouts
\usepackage{graphicx,placeins,float,amsmath,amssymb,amsthm,oubraces,cite}

\newtheorem{defn}{Definition}
\newtheorem{lem}{Lemma}

\renewcommand{\qedsymbol}{\rule{0.7em}{0.7em}}
\newenvironment{pf}[1][Lemma x]{{\it\bf Proof of #1:} }{\hfill$\qedsymbol$\\}

\def \a{ \underline{a} }
\def \b{ \underline{b} }
\def \c{ \underline{c} }
\def \d{ \underline{d} }
\def \e{ \underline{e} }
\def \f{ \underline{f} }
\def \h{ \underline{h} }
\def \u{ \underline{u} }
\def \v{ \underline{v} }
\def \w{ \underline{w} }

\def \I{\mathcal{I}}
\def \Rk { \mathbb{R}^{k}}
\def \R  { \mathbb{R} }
\def \N  { \mathbb{N} }
\def \E  { \mathbb{E} }
\def \P  { \mathbb{P} }
\def \vn { \{ \v^{n} \}_{n=0}^{\infty} }
\def \un { \{ \u^{n} \}_{n=0}^{\infty} }

\def\opt{^{\star}}

\begin{document} \onecolumn 

\title{Signal Processing Structures for Solving Conservative Constraint Satisfaction Problems}

\author{\IEEEauthorblockN{Tarek A. Lahlou and Thomas A. Baran \thanks{The authors wish to thank Analog Devices, Bose Corporation, and Texas Instruments for their support of innovative research at MIT and within the Digital Signal Processing Group.}}

\IEEEauthorblockA{Digital Signal Processing Group \\ Massachusetts Institute of Technology}}

\maketitle

\begin{abstract}
	This primary purpose of this paper is to succinctly state a number of verifiable and tractable sufficient conditions under which a particular class of conservative signal processing structures may be readily used to solve a companion class of constraint satisfaction problems using both synchronous and asynchronous implementation protocols. In particular, the mentioned class of structures is shown to have desirable convergence and robustness properties with respect to various uncertainties involving communication and processing delays. Essential ingredients to the arguments herein involve blending together functional composition methods, conservation principles, asynchronous signal processing implementation protocols, and methods of homotopy. Numerical experiments complement the theoretical presentation and connections to optimization theory are made.
\end{abstract}

\tableofcontents \newpage

\section{Introduction}
	Conservation principles in the physical sciences often play an enabling role in the predictability and tractability of system behavior at the macroscopic scale without requiring careful consideration of microscopic features at the individual subsystem level. For example, conservation of energy, a consequence of the time-invariance property of the physical principles governing an isolated system as described by Noether's theorem, allows for the analysis of physical processes without detailed handling of time-valued boundary conditions. In the context of signal processing systems, a conservative system is one for which the algorithm variables admit an organization adhering to a quadratic form of a particular class\footnote{In particular, a real symmetric bilinear form with signature of the form $(a, a, 0)$.} that is invariant with respect to the evolution of the system, e.g.~in time or space. In this way, the analysis tools of signal processing theory may naturally be utilized and extended in order to address key issues such as system stability and robustness with respect to a number of metrics for large, inhomogeneous, distributed systems that possess intrinsic conservation principles.

	Many constraint satisfaction algorithms that operate on decentralized and/or large diverse data sets are becoming increasingly burdensome for centralized and/or synchronized computing units as both the availability and geographic displacement of data continues to grow. In many contexts distributed computing environments are readily available, providing the potential for robustness with respect to fault tolerance and scalability by avoiding critical reliance upon any singular resource. The process of distributing an algorithm by organizing a non-distributed iteration to make use of such an available computational platform is generally limited in that the advantages of distribution diminish as global synchronization requirements arise, e.g.~the worst-case latency of any compute node. Another common solution is to make use of an asynchronous implementation protocol in which each node exchanges information in accordance with the computational graph without any active global organization or token passing. The issues associated with this approach may also naturally be addressed using the system analysis, organization, and implementation tools of signal processing. 
	
	The principal motivation for this paper is twofold: (1) to formally define a particular class of constraint satisfaction problems and establish their connection with a class of conservative signal processing systems and (2) to provide tractable and easily verifiable sufficient conditions for the stability of these systems in the following senses: 
	\begin{itemize}
		\item[(i)]  the systems state nears a fixed-point, i.e. a state which is invariant to system dynamics, as a consequence of deterministic and/or stochastic system evolution;
		\item[(ii)] the system returns to such an invariant state in the presence of noise and/or other perturbations.
	\end{itemize}
	In the remainder of the introduction, we formally define a conservative constraint satisfaction problem, define an associated class of conservative signal-flow structures, provide an overview of our main results, and explain how they relate to existing and recent work.

	\subsection{Conservative constraint satisfaction problems}
		Constraint satisfaction problems (CSPs) are traditionally defined as a $3$-tuple $\langle \mathcal{V}, \mathcal{D}, \mathcal{C} \rangle$ consisting of a set of variables denoted $\mathcal{V}$, a set of corresponding domains $\mathcal{D}$ over which the variables are defined, and a set of constraints between the variables denoted $\mathcal{C}$.  These may be written formally as
		\begin{eqnarray}
			\mathcal{V} &=& \left\{ v_1, \dots, v_n \right\}\\
			\mathcal{D} &=& \left\{ \mathcal{D}_1, \dots, \mathcal{D}_n \right\}\\
			\mathcal{C} &=& \left\{ \mathcal{C}_1, \dots, \mathcal{C}_r \right\},
		\end{eqnarray}
		with each variable $v_j$ satisfying $v_j\in \mathcal{D}_j$, $j=1, \dots n$, and with each $\mathcal{C}_j$, $j=1,\dots, r$ being representable as a set constraint imposed on a particular subset $\{ v_j \}$ of the variables in $\mathcal{V}$.

		\begin{defn}{\bf Conservative CSP.}\label{def:ccsp}
			We define a conservative constraint satisfaction problem (CCSP) as being reducible to a CSP described by a $3$-tuple $\langle \mathcal{V}, \mathcal{D}, \mathcal{C} \rangle$ having elements that take the following form:
			\begin{eqnarray}
				\mathcal{V} &=& \left\{ \c, \d \right\} \label{eq:ccsp1}\\
				\mathcal{D} &=& \left\{ \Rk, \Rk \right\}\\
				\mathcal{C} &=& \left\{ W, \mathcal{M} \right\},
			\end{eqnarray}
			where $W$ and $\mathcal{M}$ are constraints imposed on the entire set of variables $\left\{ \c, \d \right\}$, and where in particular $W$ is a linear subspace of $\mathbb{R}^{2k}$ that satisfies the following property:
			\begin{equation}
				\label{eq:ccspqf}
				\left\|\c\right\|^2 - \left\|\d\right\|^2 = 0,\ \ \ \left[ \begin{array}{c} \c \\ \d \end{array} \right] \in W.
			\end{equation}
		\end{defn}

		There are a variety of techniques that can be used to verify that a particular CSP is a CCSP, i.e.~that can be used to reduce a CSP to a form that satisfies \eqref{eq:ccsp1}-\eqref{eq:ccspqf}.  Among these are algebraic manipulation and reduction techniques, a key ingredient of which would be the identification and transformation of conservation principles representable as a quadratic form that is isomorphic to the left-hand side of \eqref{eq:ccspqf}, as is discussed in detail in \cite{BaranThesis}.

		The references \cite{BaranThesis,BLGS1,BLGS2,LinProg} also contain examples of a variety of practical engineering problems that are reducible to the solution of a CCSP.  These include for example the problem of solving of a broad class of linear and nonlinear optimizaion problems, as discussed in \cite{BaranThesis}.  The problem of determining the steady-state voltage and current distributions in a linear or nonlinear electrical network is also reducible to a CCSP.  There are many additional examples beyond these, a selection of which is highlighted in Section \ref{sec:examples}.

		There are generally a variety of algebraic expressions of the form of \eqref{eq:ccsp1}-\eqref{eq:ccspqf} that can be used to describe a particular CCSP.  With this in mind we focus the scope of discussion in this paper to those CCSPs for which the respective set constraints $W$ and $\mathcal{M}$ are specifically generated using functional relationships between the variables $\c$ and $\d$:
		\begin{eqnarray}
			W 			&=& \left\{ \left[ \begin{array}{c} \c    \\ G\c \end{array} \right] \ \colon \ \c \in \Rk \right\} \label{eq:define-G} \\
			\mathcal{M} &=& \left\{ \left[ \begin{array}{c} m(\d) \\ \d  \end{array} \right] \ \colon \ \d \in \Rk \right\} \label{eq:define-m}
		\end{eqnarray}
		where $G \colon \Rk \to \Rk$ is an orthogonal matrix and $m \colon \Rk \to \Rk$  is a generally nonlinear map.

	\subsection{Conservative signal processing systems}
		\begin{figure}[t!]
			\centering
		  	\centerline{\includegraphics[width=4in]{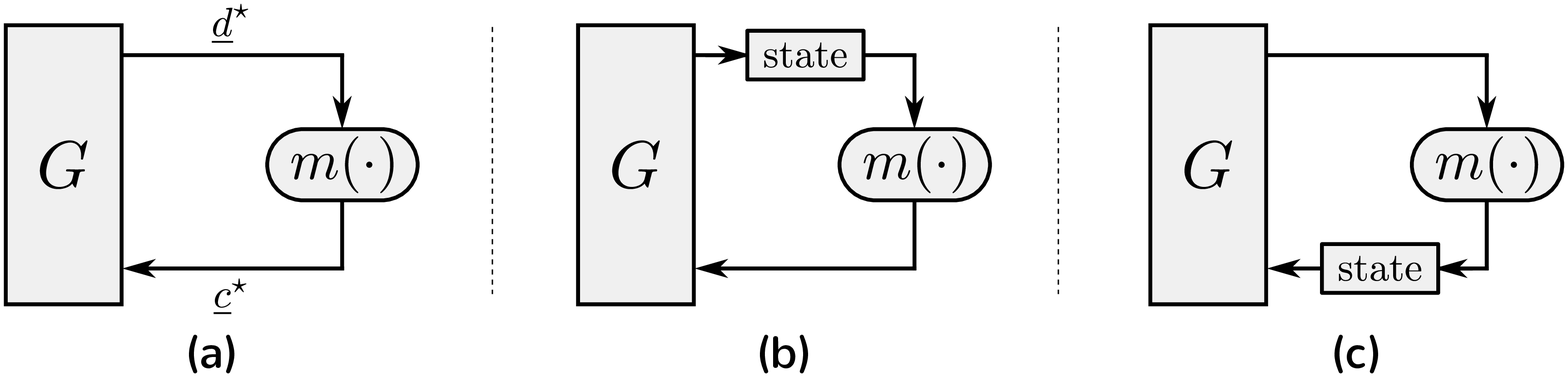}}\caption{The conservative signal-flow structure associated with (a) the algebraic system of equations \eqref{eq:c2d}-\eqref{eq:d2c} encapsulating the companion CCSP and (b)-(c) two possible realizations of (a) obtained by eliminating algebraic loops via inserting state. }
			\label{fig:conservative-signal-flow}
		\end{figure}

		From the characterization of a CCSP in Definition~\ref{def:ccsp}, it follows that an equivalent algebraic form of the problem can be posed by making use of the maps $G$ and $m(\cdot)$ in \eqref{eq:define-G} and \eqref{eq:define-m}, respectively. In particular, consider the fixed-point problem of identifying any pair of vectors $(\c\opt,\d\opt)\in\Rk\times\Rk$ which satisfy the implicit system of equations:
		\begin{eqnarray}
			\d\opt 	&	=	& 	G\c\opt 				\label{eq:c2d} \\ 
			\c\opt 	&	=	& 	m\left(\d\opt\right).	\label{eq:d2c}
		\end{eqnarray}  
		As is well-known, a signal-flow structure is, by definition, a graphical representation of a system of equations and is not necessarily indicative of a functional realization or algorithm. In the same vein, a signal-flow characterization of \eqref{eq:c2d}-\eqref{eq:d2c} as well as two possible realizations of the signal-flow structure as a computable system are depicted in Fig.~\ref{fig:conservative-signal-flow}. The utility of a conservative signal-flow system in the context of solving a CCSP is then in running the system so that the values stored in the system's state identify a fixed-point and in this sense solve the CCSP. Identifying sufficient conditions under which the synchronous and asynchronous implementation of a conservative signal-flow system does indeed unveil such a fixed-point is the primary focus of Section~\ref{sec:functional-composition}.

	\subsection{Overview of the main results}
		The following table summarizes the primary results in this paper for reference where, loosely speaking, the system operator embodying a conservative signal processing system is categorized according to Definition~\ref{def:map-types}.
		\begin{center}
			\begin{tabular}{ r | l  }
				system operator property			& reference	to sufficient condition	\\
				\hline \hline
				$\alpha$-dissipative ($\alpha < 1$) & 	Lemmas~\ref{lem:dissipative-about-fixed-point} - \ref{lem:dissipative-everywhere-filtered} \\
				passive ($\alpha = 1$) 				& 	Lemma~\ref{lem:passive-everywhere} 	\\
				$\alpha$-expansive ($\alpha > 1$)	& 	Lemma~\ref{lem:expansive-everywhere} 	\\
			\end{tabular}
		\end{center}

	\subsection{Outline of methods and organization}
		In order to present the results in this paper in such a way that avoids the critical reliance upon any particular computational platform and/or technology for realization, we restrict our treatment to identifying conditions on a function or system operator that embodies an associated signal-flow structure. As is well-known, signal-flow structures are one-to-one with generally nonlinear systems of equations and are specifically not one-to-one with any particular realization or algorithm. In this spirit, we relate the synchronous and asynchronous implementation of a signal processing system to specific methods of functional composition. The development of these relationships and their consequences with respect to establishing fixed-point properties is the primary focus of Section~\ref{sec:functional-composition} wherein a particular focus is placed on the roles of boundedness and continuity of the system operator over local and global domains. Our analysis of an asynchronous signal processing system, i.e.~a signal processing system where each delay element is an independent, randomly triggered sample-and-hold element, models the stochastic dynamics using independent Bernoulli processes. Extensions to more general discrete-time renewal processes in place of the Bernoulli processes follow in a straightforward way.

	\subsection{Manifestations of our results and relationships to recent work}
		We begin this section by qualitatively describing a limited number of manifestations of our results under various degenerations and generalizations as they appear in existing literature. It is not our intent to provide a comprehensive discussion or exhaustive list of references, but to indicate a number of interpretations of independent interest.
		\begin{itemize}
			\item Stability results for the special case of simple feedback loops which rely upon functional and relational methods, as opposed to traditional Lyapunov and storage function approaches, based on essentially generalizations of Definition~\ref{def:map-types} and its emergent properties is discussed in \cite{Zames1} and \cite{Zames2}. 

			\item The representation of systems via an interconnection of constraints has its origins in the behavioral representation theory approach to control and dynamical systems in \cite{Willems}. 
			
			\item The general approach to dissipative everywhere systems taken in this paper is akin to the use of the contraction mapping principle in the study of nonlinear equations including algebraic, integral and differential equations. Equivalently the convergence and stability results are themselves interpretable as consequences of the celebrated Banach fixed-point theorem from topology \cite{contraction}. As was previously mentioned, many of the results in this paper, including those pertaining to passive and expansive everywhere systems, may readily be extended into more general complete metric spaces over general fields using any bonafide distance metric.
			
			\item The use of operator valued random variables may be used to extend the deterministic results of the contraction mapping principle in several ways, e.g. where the operators action and/or input are stochastic in nature \cite{stochasticContraction}. Though not presented in this way, these results may correctly be used as an analysis tool for asynchronous implementations of dissipative everywhere systems. Similarly, randomization of sub-calculations, including those used in standard stochastic gradient descent methods, have been adopted as a standard tool in optimization theory \cite{hogwild}.
			
			\item The topic of weakly convergent sequences, which is discussed in standard analysis textbooks such as \cite{rudin}, may be applied in a relatively straightforward manner to a sublcass of passive everywhere system operators implemented using a synchronous implementation protocol. For direct application, constraint qualifications such as the verification of sufficient topological mixing is generally required.
			
		\end{itemize}
		
		We conclude this section by drawing attention to recent work by the present authors in which the main results of this paper find immediate application, a subset of which serve in lieu of an examples section. 
		\begin{itemize}
			\item The broader role played by conservation principles in signal processing systems, specifically as they pertain to questions surrounding identification and manipulation, is discussed at length in \cite{BaranThesis}. 
			
			\item A key connection between the class of conservative signal processing systems considered in this paper and a class of generally convex and non-convex optimization problems is established in \cite{BLGS1}. In particular, the connection relies upon conservation principles in order to invoke stationarity conditions which are closely to the principles of stationary content and co-content that can be derived using Tellegen’s theorem in electrical networks. Indeed, the key idea underlying this connection is to essentially map a conservative optimization problem into a CCSP and thus the examples in \cite{BLGS2} can readily be viewed as examples in the context of this paper.
			
			\item An in-depth case study of the use of conservative signal processing structures for synchronously and asynchronously solving linear programs with a focus on the distributed organizations and implementations is presented in \cite{LinProg}.
			
			\item Automated techniques for realizing the system operator referred to in this paper as embodying a signal processing structure is discussed in \cite{BLImpl}. In addition, an example algorithm for mapping constraints into this form is presented. 

		\end{itemize}

\newpage
\section{Synchronous and asynchronous functional composition} \label{sec:functional-composition}
	In this section various properties associated with the convergence, stability and robustness of a conservative signal processing system to synchronous and asynchronous implementation protocols are presented. The general strategy underlying the forthcoming analysis is to make use of and build upon functional composition methods in the settings of Banach and Hilbert spaces in order to provide sufficient conditions under which the sequence of system states tends to a fixed-point or invariant state. Proofs of the stated results are deferred to the appendix for brevity.

	\subsection{System operators and implementation protocols}\label{sec:composition}
		Let $T$ denote a general map from $\Rk$ into itself and let $\mathcal{F}_{T}$ denote the associated set of fixed-points, i.e.~$\mathcal{F}_{T} = \{ \v\in\Rk \colon T(\v) = \v\}$. The primary purpose of this subsection is to define the synchronous and asynchronous implementation of $T$ and operators derived from $T$ using functional composition methods in anticipation of the fixed-point results of Subsections~\ref{sec:dissipative} through \ref{sec:general}. To this end, we proceed with the implicit assumption that $T$ is a {\it system operator} and thus embodies a signal processing system by which we mean that $T$ maps a vector $\v$ comprised of the current values in the signal processing systems state to the vector $T(\v)$ corresponding to the action of the constraints imposed by the signal processing system on those values at the input to the systems state. This assumption is taken for contextual convenience in order to focus the terminology introduced next and is not necessary for technical reasons. Furthermore, we note that $\mathcal{F}_{T}$ may generally be empty, singleton, convex, the union of several disjoint subsets of $\Rk$, etc., depending on any number of properties belonging to the associated system operator $T$. We shall consequently make clear when we have explicitly made an assumption about $\mathcal{F}_{T}$ in order to restrict the class of system operators considered as opposed to when a property of the class of system operators at hand implies a property of $\mathcal{F}_{T}$.

		We begin by defining a homotopy map for use with a continuation scheme wherein the general idea is to continuously deform a simple system operator with desirable convergence properties into the given system operator $T$. The utility of such a scheme in the context of fixed-point analysis is in sequentially solving or approximately solving the family of deformed system operators for a fixed-point in such a manner that the final solution obtained is a fixed-point of the original operator. The function which describes these deformations is referred to as the homotopy map of the system. In this paper we restrict ourselves to only one analytic form of a homotopy map which appears in the following definition although alternative parameterizations may generally be valid homotopy maps as well.

		\begin{defn}{\bf Homotopic continuation system operator \boldmath$T_{h}$} \label{def:homotopy-system} \\ 
			Let $T_{(1)} \colon \Rk \to \Rk$ denote a system operator and let $T_{(0)} \colon \Rk \to \Rk$ denote a second system operator for which the homotopy map or homotopic continuation system operator $T_{h} \colon [0,1] \times \Rk \to \Rk$ for the system operator $T_{(1)}$ is:
			\begin{eqnarray}
				T_{h}\left( \rho, \v \right) \triangleq \rho T_{(1)}(\v) + (1-\rho)T_{(0)}(\v). \label{eq:homotopy-system}
			\end{eqnarray}
		\end{defn}
		
		The family of fixed-point problems defined by the homotopy map $T_{h}$ is given by $T_{h}(\rho,\v)=\v$ for $\rho\in[0,1]$ and the standard implementation strategy is to track the fixed-points $T_{h}(\rho,\v)=\v$ starting from $(\rho,\v)=(0,\v^*)$ as $\rho$ progresses from $0$ to $1$ where $v^*$ is an easily obtained fixed-point of $T_{(0)}$. Typical realizations of this strategy discretize $\rho$ and sequentially increment it by some sufficiently small amount where, under appropriate regularity conditions, this procedure will result in a pair $(\rho,\v)=(1,\v^\star)$ where $T_{(1)}(\v^\star)=\v^\star$ and thus any fixed-point of $T_{h}(1,\cdot)$ is also a fixed-point of $T_{(1)}$. Observe that the trivial selection of the second system operator $T_{(0)}$ as the identity operator on $\Rk$ results in a special case of the homotopic continuation system operator that is able to identify $\v^\star \in \mathcal{F}_{T_{(1)}}$  after identifying a fixed-point of $T_{h}$ for any non-zero value of $\rho$. So, in light of this and the observation that $\rho$ need not be contained to the unit interval for this to be true, we define next a filtered system operator which subsumes these observations.

		\begin{defn}{\bf Filtered system operator \boldmath$T_{f}$} \label{def:filtered-system} \\
			Let $T \colon \Rk \to \Rk$ denote a system operator from which the filtered system operator $T_{f} \colon (0,\infty) \times \Rk \to \Rk$ is:
			\begin{equation} \label{eq:filtered-system}
				T_{f}\left( \rho, \v \right) \triangleq \rho T ( \v ) + ( 1 - \rho ) \v. 
			\end{equation}
		\end{defn}
		
		We now justify the observation above. Denote by $\v^\star \in \Rk$ an invariant state of $T_f$ for a fixed value of $\rho$, i.e.~$\v^\star \in \mathcal{F}_{T_f}$. After some straightforward manipulations we obtain that $T( \v^\star ) = \v^\star$ and thus it follows that any fixed-point of $T_f$ is also a fixed-point of $T$ since the analysis holds for any selection of $\rho > 0$. We note for completeness that the converse of this fact is also true, but is not of use in this paper. One signal processing interpretation of a filtered system operator is that the system's state is updated with a weighted average or convex combination (if $\rho \in (0,1)$) of the values currently in the system's state and the action of $T$ on those values. This procedure is a common weight adjustment protocol in adaptive filtering systems. 

		Prior to defining the implementation of a signal processing system via synchronous and asynchronous functional composition, we first comment on the generality which may be assumed by our use of the system operator $T$. In particular, without loss of generality, $T$ may itself be the composition of finitely many system operators, i.e.~$T$ may correspond to a complete iteration of an iterated function system (without the usual restriction of each map being a contraction). Written more formally, for a signal processing system decomposed into the sequential iterated function system described by:
		\begin{eqnarray} \label{eq:iterated-function-system}
			&\left\{ T_{(i)}: \Rk \to \Rk \colon i = 1, 2, \dots, m\right\}, & m \in \mathbb{N},
		\end{eqnarray}
		where the implementation of each map $T_{(i)}$ may make use of either of the system operators previously discussed, the system operator $T$ then corresponds to the composite map $T = T_{(m)} \circ T_{(m-1)} \circ \cdots \circ T_{(1)}$. Subsequently, the following definitions are made explicitly using the system operator $T$ with an understanding that the specific form and complexity of $T$, e.g. a filtered system operator $T_f$, may be taken in place of $T$ in the actual realization of the system. 

		\begin{defn}{\bf Synchronous implementation protocols} \label{def:synchronous} \\
			Let $T \colon \Rk \to \Rk$ denote a system operator. The state evolution sequence associated with a synchronous implementation of $T$ starting from some initial system state $\v^{0} \in \Rk$ is defined as the sequence $\vn$ generated according to:
			\begin{eqnarray} \label{eq:synchronous}
				& \v^{n} \triangleq T(\v^{n-1}), & n \in \N. 
			\end{eqnarray}
		\end{defn}
		
		From the intuitive notion of an asynchronous signal processing system as one for which the system state acts, for each scalar delay or state element independently, as a randomly triggered sample-and-hold element, we next define an asynchronous implementation protocol as the discrete-time model of this behavior where each delay element is driven by an independent discrete-time Bernoulli processes.  

		\begin{defn}{\bf Asynchronous implementation protocols} \label{def:asynchronous} \\
			Let $T \colon \Rk \to \Rk$ denote a system operator. The state evolution sequence associated with an asynchronous implementation of $T$ starting from some initial, deterministic system state $\v^{0} \in \Rk$ is defined as the sequence $\vn$ generated according to:
			\begin{eqnarray} \label{eq:asynchronous}
				& \v^{n} \triangleq D^{(p)}T(\v^{n-1}) + (I_k-D^{(p)})\v^{n-1}, & n \in \N, 
			\end{eqnarray}
			where $I_k$ is the $k \times k$ identity matrix and $D^{(p)}$ is a $k \times k$ stochastic, binary, diagonal matrix whose $k$ diagonal elements are i.i.d., Bernoulli and independent of $\v^{n-1}$, taking values $D^{(p)}_{ii} = 1$ with probability $p$ and $D^{(p)}_{ii} = 0$ with probability $1-p$.
		\end{defn}

		Observe that the update mechanism used to generate the state evolution sequence for an asynchronous implementation of the system operator $T$ closely resembles that of the synchronous implementation of the filtered system operator $T_f$. In particular, modifying the behavior of the stochastic matrix $D^{(p)}$ to its expectation yields precisely the deterministic state evolution update procedure of the filtered system operator where the filter coefficient $\rho$ takes the value of the probability of an asynchronous delay element firing for each $n$. By extension, this procedure is also closely related to the synchronous implementation of the homotopic continuation system operator where $T_{(0)} = I_k$.

		Having established two implementation protocols for which the dynamics of the associated signal processing system are governed by deterministic and stochastic state elements, we now make the meaning of a convergent state evolution sequence precise in each setting as well as establish an associated notational convention. In particular, we shall write that the state evolution sequence $\vn$ converges to an invariant state $\v^\star$ using the notation:
		\begin{eqnarray} \label{eq:convergence}
			& \v^n \xrightarrow{(\ell,r)} \v^\star &
		\end{eqnarray}
		where the meaning of convergence is specific to the implementation protocol and is described by the parameters $(\ell,r)$  defined shortly. For state evolution sequences produced using Definition~\ref{def:synchronous}, convergence is in the sense of a metric space equipped with the Euclidean distance metric. However, in order to establish either a convergence rate or order depending on the result obtained, we parameterize \eqref{eq:convergence} using $\ell=1$ or $2$. More formally, \eqref{eq:convergence} in the synchronous setting precisely means that:
		\begin{equation*}
		 	\text{for any $\epsilon > 0$ there exists an $n_0 \in \N$ such that $\|\v^n - \v^\star\|^\ell < \epsilon$ for all $n > n_0$}
		\end{equation*}
		which we rewrite succinctly as:
		\begin{eqnarray} \label{eq:convergence-synch}
			\lim_{n \to \infty} \left\| \v^n - \v^\star \right\|^\ell = 0. 
		\end{eqnarray}
		Similarly, for state evolution sequences produced using Definition~\ref{def:asynchronous}, convergence is in the sense of $r$-th mean in a metric space for $r = 1$ or $2$. More formally,  \eqref{eq:convergence} in the asynchronous setting precisely mean that:
		\begin{equation*}
			\text{for any $\epsilon > 0$ there exists an $n_0 \in \N$ such that $\E\left[ \|\v^n - \v^\star\|^{r}\right] < \epsilon$ for all $n > n_0$}
		\end{equation*}
		which we rewrite succinctly as:
		\begin{eqnarray} \label{eq:convergence-asynch}
			\lim_{n \to \infty} \E \left[ \left\| \v^n - \v^\star \right\|^{r} \right] = 0
		\end{eqnarray}
		where $\E[\cdot]$ is the expectation operator. Therefore, in the synchronous setting, we say that $\vn$ convergnes linearly to $\v^\star$ at rate $\mu\in(0,1)$ and that $\vn$ converges to $\v^\star$ at order $p$ provided that:
		\begin{eqnarray}
			\lim_{n\to\infty}\frac{\|\v^n - \v^\star\|}{\|\v^{n-1}-\v^\star\|} = \mu, & \text{ and } & \|\v^n - \v^\star\| < \frac{K}{n^p},\,\,\,\,\,n\geq0,
		\end{eqnarray}
		respectively, where $K$ is an arbitrary constant. We use an analagous definition for the asynchronous setting. Moreover,  convergence in mean square $(r=2)$ has many implications which are useful in understanding the dynamics of an asynchronous signal processing system \cite{RVConvergence}. We now state a few such implications for completeness:
		\begin{itemize}
			\item[(i)] 		convergence in mean square ($r=2$) implies convergence of a subsequence of $\vn$ almost surely;
			\item[(ii)] 	convergence in mean square ($r=2$) implies convergence in mean ($r=1$) by application of Jensen's inequality;
			\item[(iii)]	convergence in mean ($r=1$) implies convergence in probability by application of Markov's inequality;
			\item[(iv)]		convergence in probability implies convergence in distribution.
		\end{itemize}
		Convergence of $\vn$ in mean square ($r=2$) to a deterministic state $\v^\star$ implies that the variance of $\v^n$ tends to zero (by application of the Cauchy-Schwarz inequality). Note that convergence in mean square does not necessarily preclude convergence of $\vn$ almost surely.

		In anticipation of the convergence analysis for asynchronous implementations, we now justify a useful identity which we shall freely utilize in our proofs. Specifically:
		\begin{eqnarray}
			\E\left[\|\v^n - \v^\star \|^{r}\right] & = & p \E\left[\left\|\v^n - \v^\star \right\|^r \mid D^{(p)} = I_k  \right] + (1-p)\E\left[\left\|\v^n - \v^\star \right\|^r \mid D^{(p)} = 0  \right]  \label{eq:asynchronous-identity-pre} \\
			& = & p\E\left[\left\|T(\v^{n-1}) - T(\v^\star) \right\|^r  \right] + (1-p)\E\left[\left\|\v^{n-1} - \v^\star \right\|^r  \right] \label{eq:asynchronous-identity}
		\end{eqnarray}
		where the first equality is due to the law of total expectation and the second equality is due to \eqref{eq:asynchronous} and the definition of a fixed-point. Note that \eqref{eq:asynchronous-identity} is only valid for a direct implementation of $T$ while \eqref{eq:asynchronous-identity-pre} may be used to generate the associated expressions for $T_f$ and $T_h$ and does not rely upon $\v^\star$ being a fixed-point.

		Let $B(\c, r) \subset \Rk$ describe the closed Euclidean ball of radius $r > 0$ and center $\c \in \Rk$ given by:
		\begin{equation} \label{eq:closed-bowl}
			B(\c, r) \triangleq \left\{ \v \in \Rk \colon \left\| \c - \v \right\| \leq r \right\}. 
		\end{equation}
		If a state evolution sequence generated by a synchronous or asynchronous implementation of a system operator satisfies:
		\begin{eqnarray} \label{eq:inside-ball}
			&\|\v^n - \c\|\leq r \hspace{2em}\text{ or }\hspace{2em} \E\left[\|\v^n - \c\|\right] \leq r, &  n = 0, 1, 2, \dots,
		\end{eqnarray} 
		respectively, then we shall say that $\vn$ is contained within $B(\c,r)$.

\begin{figure}[t]
			\centering
		  	\centerline{\includegraphics[width=6.25in]{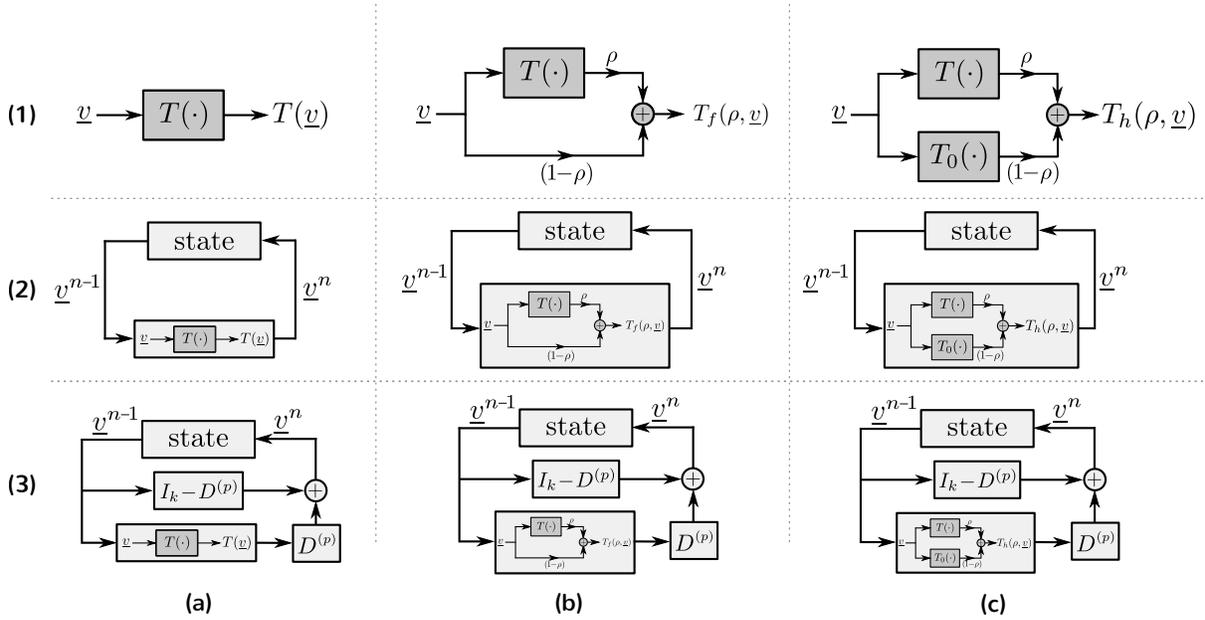}}\caption{An illustration of the system operator and implementation protocol definitions in Section~\ref{sec:composition}. In particular, column (a) pertains to a system operator $T$, column (b) pertains to the filtered system operator $T_{f}$ in Definition~\ref{def:filtered-system} for a fixed-value of $\rho$, column (c) pertains to the homotopic continuation system operator $T_{h}$ in Definition~\ref{def:homotopy-system} where $\rho$ may be adaptively selected during runtime, row (1) depicts a signal-flow representation of the associated system operators, row (2) illustrates a synchronous implementation of each system operator as described by Definition~\ref{def:synchronous}, and row (3) illustrates an asynchronous implementation of each system operator as described by Definition~\ref{def:asynchronous} where the depicted system state is that of traditional vector delay elements and the stochastic operator  $D^{(p)}$ enforces the asynchronicity.}\vspace{-0.0625in}
			\label{fig:composition-definitions}
		\end{figure}

		With these definitions in place, we remark that the actual realization of a given system operator $T$ in practice using either implementation protocol may often only be mathematically characterized by \eqref{eq:synchronous} or \eqref{eq:asynchronous}. As a straightforward example, when $T$ represents a general nonlinear coordinatewise map the corresponding asynchronous implementation may only compute the randomly selected subset of coordinates described by $D^{(p)}$ for each $n$ rather than computing all $k$ and discarding the complement of that subset, as is suggested by \eqref{eq:asynchronous}. The definitions in this subsection are summarized in Figure~\ref{fig:composition-definitions} using signal-flow representations. 
		
	\subsection{A categorization of system operators}
		In the sequel, we shall make frequent use of the categories of system operators defined next in order to assign convergence and stability properties to various implementations of a system operator $T$. These definitions, and the analysis which follows, make explicit use of finite-dimensional real linear spaces equipped with the standard Euclidian distance metric, as this setting is sufficient for our purposes and is consistent with accepted notions of convergence in numerical analysis. Nonetheless, we note from the outset that the presented results may be translated into more general settings in a straightforward way.

		\begin{defn}{\bf {\boldmath$\alpha$}-conic system operators} \label{def:alpha-conic} \\
			A system operator $T \colon \Rk \to \Rk$ is called $\alpha$-conic about $\v \in \Rk$ provided that there exists a real constant $\alpha \geq 0$ such that: 
			\begin{eqnarray} \label{eq:alpha-conic}
				\sup_{\u \neq \underline{0}} \frac{\left\|T(\v+\u) - T(\v) \right\|}{\left\| \u \right\| } \leq \alpha. 
			\end{eqnarray}
			Furthermore, the system operator is called $\alpha$-conic everywhere if it is $\alpha$-conic about all elements $\v \in \Rk$.
		\end{defn}

		We specifically draw attention to the relationship between \eqref{eq:alpha-conic} and continuity. In particular, the definition of a system operator which is $\alpha$-conic everywhere is equivalent to the standard condition of Lipschitz continuity with Lipschitz constant $\alpha$ while a system operator satisfying only $\alpha$-conicity around some element $\v$ is weaker since bounded discontinuities away from $\v$ are admissible.  Specifically, for an $\alpha$-conic everywhere system operator $T$, observe that the inequality $\left\| T(\v + \u) - T(\v) \right\| \leq \alpha \left\| \u \right\|$ is implied by \eqref{eq:alpha-conic} since by assumption $T$ is $\alpha$-conic for all $\v \in \Rk$ and $\alpha$ is the least upper bound for all $\u \neq \underline{0}$. Re-writing this inequality using the translated arguments $\v'$ and $\u'$ where $\v' = \v + \u$ and $\u'= \v$ yields:
		\begin{eqnarray} \label{eq:alpha-conic-everywhere-equivalent} 
			 \left\| T(\v') - T(\u') \right \| \leq \alpha \left\| \v' - \u' \right\|, & & \forall \v', \u' \in \Rk,
		\end{eqnarray}
		hence \eqref{eq:alpha-conic-everywhere-equivalent} is an equivalent condition for a system operator $T$ to be $\alpha$-conic everywhere with the same parameter $\alpha$.

		We shall frequently postulate a system operator $T$ which is $\alpha$-conic about some element $\v \neq \underline{0}$ for which it may be helpful from either an intuitive or analytic perspective to consider a related system operator which is $\alpha$-conic about $\underline{0}$. The following lemma permits us to do just this, i.e.~to extend arguments related to convergence and stability pertaining to the class of $\alpha$-conic about $\underline{0}$ system operators to those which are $\alpha$-conic about arbitrary $\v$ and vice versa without the loss of any generality. 

		\begin{lem} \label{lem:translate-to-zero}
			Let $T_{(1)}$ denote a system operator which is $\alpha$-conic about $\v^1$ with an associated non-empty set of fixed-points $\mathcal{F}_{T_{(1)}} = \{ \v \in \Rk \colon T_{(1)}(\v) = \v\}$. Then there exists a related system operator $T_{(2)}$ which is $\alpha$-conic about $\underline{0}$ with an associated non-empty set of fixed-points $\mathcal{F}_{T_{(2)}} = \{ \v \in \Rk \colon T_{(2)}(\v) = \v \}$ such that an analytic map exists relating the sets of fixed-points $\mathcal{F}_{T_{(1)}}$ and $\mathcal{F}_{T_{(2)}}$ in such a way that identifying a fixed-point of either system identifies a fixed-point of both.
		\end{lem}	
		For the specific system operator utilized in the proof in \eqref{eq-lem:translate-to-zero-10}, it follows that $T_{(2)}(\underline{0})=\underline{0}$ if and only if $\v^1 \in \mathcal{F}_{T_{(1)}}$. 

		Next, we partition the class of $\alpha$-conic system operators introduced in Definition~\ref{def:alpha-conic} into three subclasses which are, roughly speaking, used to assign the convergence and stability properties discussed in Sections~\ref{sec:dissipative} through \ref{sec:general}.

		 \begin{defn}{\bf  \boldmath{$\alpha$}-dissipative, passive, and  \boldmath{$\alpha$}-expansive system operators} \label{def:map-types} \\
			Let $T \colon \Rk \to \Rk$ denote a system operator which is $\alpha$-conic about some element $\v \in \Rk$. We shall refer to $T$ as being:
			\begin{itemize}
				\item[(i)] 		$\alpha$-dissipative about $\v$ provided that $\alpha \in [0, 1)$;  

				\item[(ii)]		passive about $\v$ provided that $\alpha = 1$; 

				\item[(iii)] 	$\alpha$-expansive about $\v$ provided that $\alpha > 1$. 
			\end{itemize}
			Furthermore, when $T$ is $\alpha$-conic everywhere, i.e.~satisfies \eqref{eq:alpha-conic} for all $\v \in \mathbb{R}^k$, we shall refer to $T$ as being:
			\begin{itemize}
				\item[(i)] 		$\alpha$-dissipative everywhere provided that $\alpha \in [0, 1)$;  

				\item[(ii)]		passive everywhere provided that $\alpha = 1$; 

				\item[(iii)] 	$\alpha$-expansive everywhere provided that $\alpha > 1$. 
			\end{itemize}
		\end{defn}

		In the terminology of Lipschitz continuous functions, any system operator which is $\alpha$-dissipative or passive everywhere is also referred to as being contractive or non-expansive, respectively.

		We now discuss the invariance of the $\alpha$-conic everywhere property with respect to composition by considering an iterated function system of the form \eqref{eq:iterated-function-system}. The utility of this invariance, formalized in the following lemma, is in potentially obtaining an $\alpha$-dissipative everywhere system operator from the composition of a single dissipative everywhere system operator with any number of passive and/or expansive everywhere system operators. 

		\begin{lem} \label{lem:conic-composition}
			Let $\left\{ T_{(i)}: \Rk \to \Rk \colon i = 1, 2, \dots, m\right\}$ be a set of $m\in\N$ $\alpha_{i}$-conic everywhere system operators. The system operator $T_{m} \circ T_{m-1} \circ \cdots \circ T_{1}$ is $\alpha$-conic everywhere with parameter $\displaystyle \prod_{i=1}^{m} \alpha_i$.
		\end{lem}

		We next define the class of system operators which act as isometries on $\mathbb{R}^{k}$.

		\begin{defn}{\bf Neutral system operators} \label{def:neutral} \\ 
			A system operator $T \colon \Rk \to \Rk$ is said to be neutral provided that:
			\begin{eqnarray} \label{eq:neutral}
				&\left\| T(\v) \right\| = \left\| \v \right\|,  & \forall \, \v \in \Rk. 
			\end{eqnarray}		
		\end{defn}

		Any map which is known to be linear in addition to being neutral is necessarily an orthogonal transformation, an important subset of linear isometries. Finally, we define a subset of affine system operators which may be interpreted as introducing an additive bias into the associated signal processing system. 

		\begin{defn}{\bf Source system operators} \label{def:source} \\
			A system operator $T \colon \Rk \to \Rk$ is said to be a source provided that it has the algebraic form:
			\begin{eqnarray} \label{eq:source}
				T(\v) = S\v + \e
			\end{eqnarray}
			 where $S$ is a linear neutral map on $\Rk$ and $\e \in \Rk$ is a constant.
		\end{defn}

		Neutral and source system operators are readily shown to be passive everywhere. Therefore, by application of Lemma~\ref{lem:conic-composition} it follows that if any number of neutral and source system operators are introduced into an iterated function system corresponding to an overall $\alpha$-dissipative everywhere map then the new overall system operator remains $\alpha$-dissipative everywhere with the same parameter.

	\subsection{Convergence and stability results related to $\alpha$-dissipative system operators} \label{sec:dissipative}
		In this section a number of sufficient conditions under which both synchronous and asynchronous implementations of various classes of system operators associated with $\alpha$-dissipativity in some manner are established. Loosely speaking, the system operator properties considered in this section are generally sufficient to ensure that either a unique fixed-point exists or that when a fixed-point is assumed to exist it is unique. We specifically call attention to the fact that fixed-point properties belong to the system operator itself and are independent of the implementation protocol details used in establishing them.

		Consider the class of system operators which are $\alpha$-dissipative about some state $\v$. This restriction does not require any continuity of the system operator except at the given state $\v$ meaning that bounded discontinuities away from $\v$ are admissible. Further, we note that being $\alpha$-dissipative at $\v$ does not preclude the system operator from being $\alpha$-expansive about some other states. The following two lemmas consider this class of system operators when $\v$ is taken to be and not be a fixed-point, respectively. Indeed, if a system operator is known to be $\alpha$-dissipative about a fixed-point it follows that the associated signal processing system embodied by the system operator will settle linearly to an invariant state regardless of its initial configuration.  While the existence of a fixed-point of the system operator is assumed for this scenario, its uniqueness follows as a consequence of \eqref{eq:alpha-conic}. We summarize these remarks into the following lemma. 
		
		\begin{lem} \label{lem:dissipative-about-fixed-point}
			Let $T$ denote a system operator which is $\alpha$-dissipative about $v^\star$ where $\v^\star$ is a fixed-point of $T$ and let $\v^0$ be an arbitrarily selected, deterministic initial system state. Then the state evolution sequence $\vn$ generated by either a synchronous or asynchronous implementation of $T$ satisfies $\v^n \xrightarrow{(1,2)} \v^\star$. Moreover, $\v^\star$ is the unique fixed-point of $T$.
		\end{lem}

		It is typically difficult to verify Lemma~\ref{lem:dissipative-about-fixed-point} in practice since it involves both the assumption that a fixed-point exists as well as knowledge of the behavior of the system operator around it. Alternatively, it is more common that a system operator is known to be $\alpha$-dissipative about some point $\c\in\Rk$ which is not a fixed-point.  Motivated by this, the next lemma ensures that the state evolution sequence associated with such a system operator will be contained within a closed Euclidean ball centered around $\c$ after finite transient effects die off. 

		\begin{lem} \label{lem:dissipative-about-some-point}
			Let $T$ denote a system operator which is $\alpha$-dissipative about $\c$ and let $\v^0$ be an arbitrarily selected, deterministic initial system state. Then, for every $\epsilon > 0$, the state evolution sequences generated by synchronous and asynchronous implementations of $T$ are, after a finite number of initial iterations, forever contained within the closed Euclidean ball $B\left(\c, \frac{\left\| T(\c) - \c \right\| }{1-\alpha} + \epsilon\right)$.
		\end{lem}

		The entrapment result stated above neither assumes nor concludes anything about a fixed-point of the system operator $T$ since in general a fixed-point need not exist if $T$ is only known to be $\alpha$-dissipative about an arbitrary state $\c$. However, if we take $\c$ to be a fixed-point we find that Lemma~\ref{lem:dissipative-about-some-point} is consistent with Lemma~\ref{lem:dissipative-about-fixed-point} in the sense that the radii achieved for the synchronous and asynchronous implementations is given by $\epsilon$ and may be selected to be arbitrarily small. This type of agreement is typical of finite-time convergence results. In Section~\ref{sec:general} the difference in the utility of  Lemma~\ref{lem:dissipative-about-fixed-point} and Lemma~\ref{lem:dissipative-about-some-point} in assembling convergence results for general system operators is elaborated on in the context of an example.
	
		We now consider the class of $\alpha$-dissipative everywhere system operators and remark that the consequences are considerably strengthened. For example, in comparison with Lemma~\ref{lem:dissipative-about-fixed-point}, a fixed-point of an $\alpha$-dissipative everywhere system operator need not be assumed while similar convergence results and the uniqueness of the fixed-point are still guaranteed. The signal processing system associated with the $\alpha$-dissipative everywhere system operator will settle to a unique invariant state for any initial configuration and this convergence is linear for both implementation protocols and in particular with rate $\alpha$ in the synchronous setting. The following lemma subsumes these remarks. 

		\begin{lem}\label{lem:dissipative-everywhere}
			Let $T$ denote a system operator which is $\alpha$-dissipative everywhere and let $\v^0$ be an arbitrarily selected, deterministic initial system state. Then the state evolution sequence $\vn$ generated by either a synchronous or asynchronous implementation of $T$ satisfies $\v^n \xrightarrow{(1,2)} \v^\star$ where $\v^\star\in\mathcal{F}_{T}$. Moreover, $\v^\star$ is the unique fixed-point of $T$.
		\end{lem}
		
		Next, we treat the class of system operators which form the middle ground between those considered in Lemmas~\ref{lem:dissipative-about-fixed-point}-\ref{lem:dissipative-about-some-point} ($\alpha$-dissipative about a point) and Lemma~\ref{lem:dissipative-everywhere} ($\alpha$-dissipative everywhere). In particular, we consider system operators which are $\alpha$-dissipative over a Euclidean ball in $\Rk$ but not necessarily everywhere. In this case, if the system operator is properly initialized with a state from inside the ball and additionally satisfies a distance preserving property involving the translation of the center of the ball by the system operator then convergence to an invariant state is still guaranteed. We restate this in the following lemma and proceed with a closed ball in our presentation with an understanding that the uniform continuity of $T$ over an open ball inherited from \eqref{eq:alpha-conic} is sufficient to extend it to the closure of the ball while retaining the same parameter $\alpha$.

		\begin{lem} \label{lem:dissipative-everywhere-ball}
			Let $T$ denote a system operator which is $\alpha$-dissipative about all $\v \in B(\c, r)$. Then the state evolution sequence $\vn$ generated by either a synchronous or asynchronous implementation of $T$ satisfies $\v^n \xrightarrow{(1,1)} \v^\star$ where $\v^\star$ is the unique fixed-point of $T$ in $B(\c,r)$ provided that $\v^0 \in B(\c,r)$ and: 
			\begin{equation}\label{eq:dissipative-everywhere-ball-condition}
				\left\| \c - T(\c) \right\| \leq (1-\alpha)r.
			\end{equation}
		\end{lem}

		This result is consistent with Lemma~\ref{lem:dissipative-everywhere} in the sense that $\Rk$ can be thought of as a Euclidean ball with an arbitrarily large radius. In this sense, the bound in \eqref{eq:dissipative-everywhere-ball-condition} can also be made arbitrarily large hence the location of the initial system state $\v^0$ becomes insignificant. Certifying that a given system operator satisfies Lemma~\ref{lem:dissipative-everywhere-ball} has the convenient property that it only requires a single application of $T$. Furthermore, by application of Lemma~\ref{lem:translate-to-zero}, a related system operator may be defined for which the condition \eqref{eq:dissipative-everywhere-ball-condition} is specifically that the operator norm of the related system operator must be upper bounded by $r(1-\alpha)$.

		Beyond the proof of Lemma~\ref{lem:dissipative-everywhere-ball}, the application of Lemma~\ref{lem:dissipative-everywhere} is useful in a number of additional contexts. One such case, formalized in the next lemma, pertains to general system operators which themselves are not necessarily $\alpha$-conic but for which some finite self-composition is $\alpha$-dissipative everywhere. Indeed, the following lemma states that this condition is sufficient for both synchronous and asynchrnous implementations of $T$ to behave as if $T$ satisfied Lemma~\ref{lem:dissipative-everywhere}.

		\begin{lem}\label{lem:dissipative-everywhere-composition}
			Let $T$ denote a system operator and let $\v^0$ be an arbitrarily selected, deterministic initial system state. If there exists a finite $m \in \N$ for which $T^m = T \circ T^{m-1}$ is $\alpha$-dissipative everywhere then the state evolution sequence $\vn$ generated by either a synchronous or asynchronous implementation of $T$ satisfies $\v^n \xrightarrow{(1,2)} \v^\star$ where $\v^\star$ is the unique fixed-point of $T$. 
		\end{lem}

		The provision that a composite system operator $T^m$ be $\alpha$-conic everywhere does not imply that $T$ is, nor does it imply that $T$ is continuous in any sense. As a straightforward example of the former, consider the system operator $T:\R \to \R$ given by $T(v) = e^{-v}$ which is readily verified to be a passive system operator over the domain of positive reals, however $T$ does not satisfy \eqref{eq:alpha-conic} everywhere for any finite value of $\alpha$. On the other hand, the composite map $T \circ T(v) = e^{-e^{-v}}$ does with $\alpha = \frac{1}{e}$.

		The classes of system operators associated with $\alpha$-dissipativity considered in this section have not needed to utilize the modified system operators in Definitions~\ref{def:homotopy-system} and \ref{def:filtered-system} in order to establish convergence guarantees. In general, the utility of a filtered system operator $T_f$ with parameter $\rho \in (0,1)$ will be made clear in the coming sections, however a key benefit may still be found in the dissipative setting with respect to accelerating convergence by taking $\rho > 1$. The next lemma presentes this, i.e.~that convergence in the sense of Lemma~\ref{lem:dissipative-everywhere} may be achieved at a better rate by proper selection of $\rho$ where $\rho$ plays a role akin to a fixed step length selection rule of an iteration scheme.

		\begin{lem} \label{lem:dissipative-everywhere-filtered}
			Let $T$ denote a system operator which is $\alpha$-dissipative everywhere and let $\v^0$ be an arbitrarily selected, deterministic initial system state. Then the state evolution sequence $\vn$ generated by either a synchronous or asynchronous implementation of the filtered system operator $T_f$ satisfies $\v^n \xrightarrow{(1,2)} \v^\star$ where $\v^\star\in\mathcal{F}_{T}$ is unique provided that $\rho \in \left(0, \frac{2}{1+\alpha} \right)$.
		\end{lem}

	\begin{figure}[t!]
		\centering
	  	\centerline{\includegraphics[width=\textwidth]{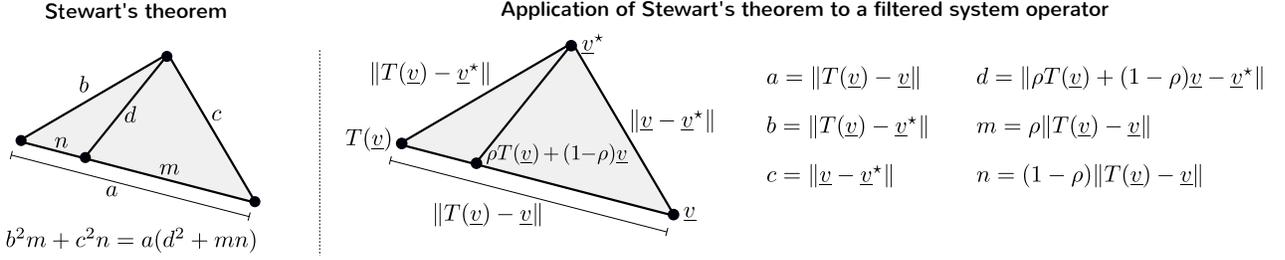}}\caption{An illustration of the application of Stewart's theorem (on the left) to the update procedure associated with a filtered system operator $T_f$ (on the right) where $\v^{n-1} = \v$ and $\v^{n} = \rho T(\v) + (1-\rho)\v$. The depicted sizes $b$ and $c$ are equal, consistent with the the passive everywhere property of $T$.}
		\label{fig:passive-everywhere-triangle}
	\end{figure}

	\subsection{Convergence and stability results related to passive system operators} \label{sec:passive}
		In this section we consider the class of passive everywhere system operators. Unlike the $\alpha$-dissipative system operators considered in Section~\ref{sec:dissipative} wherein a unique fixed-point is typically guaranteed to exist, a fixed-point of a system operator satisfying Definition~\eqref{def:alpha-conic} everywhere for $\alpha = 1$ may or may not exist, and when one does exist it may or may not be unique. However, when $\mathcal{F}_T$ is non-empty, it is a convex set, i.e.~if $\u^\star$ and $\v^\star$ are two different fixed-points of $T$ then so is the continuum of states given by $\gamma\v^\star + (1-\gamma)\u^\star$ for $\gamma \in[0,1]$. To see this, let $T$ denote a passive system operator with fixed-points $\u^\star$ and $\v^\star$ and observe that the inequality: 
		\begin{eqnarray}
			\left\| \u^\star - \v^\star \right\| & \leq & \left\| \u^\star - T\left(\gamma\v^\star + (1-\gamma)\u^\star\right) \right\| + \left\| T\left(\gamma\v^\star + (1-\gamma)\u^\star\right) - \v^\star \right\| \\
			& \leq & \left\| \gamma \left( \u^\star - \v^\star\right) \right\| + \left\|\left( 1 - \gamma \right) \left( \u^\star - \v^\star \right) \right\| \\
			& \leq & \left\| \u^\star - \v^\star \right\|
		\end{eqnarray}
		actually holds with equality where $\gamma \in [0,1]$ and the first inequality is due to the subadditivity of norms, the second inequality is due to the passivity of $T$, and the third inequality is due to the homogeneity of norms. The implication of the overall equality after application of the triangle inequality is that $T(\gamma\v^\star + (1-\gamma)\u^\star)$ and $\gamma\v^\star + (1-\gamma)\u^\star$ lie on the line segment connecting $\v^\star$ and $\u^\star$  and subsequently that $\gamma\v^\star + (1-\gamma)\u^\star \in \mathcal{F}_{T}$. Hence, we have shown that $\mathcal{F}_{T}$ is a convex set. When $\mathcal{F}_{T}$ is singleton this is trivially true. Notice the minimal assumption required for this argument is that $T$ be passive about $\mathcal{F}_{T}$ rather than passive everywhere.

		When a passive everywhere system operator possesses a non-empty convex set of fixed-points, it follows that the associated signal processing system will settle to an invariant state regardless of its initialization so long as the system is implemented using a filtered system operator with parameter $\rho$ in the open unit interval. The following lemma formalizes this fact.
		
		\begin{lem} \label{lem:passive-everywhere}
			Let $T$ denote a system operator which is passive everywhere with a non-empty convex set of fixed-points $\mathcal{F}_{T}$ and let $\v^0$ be an arbitrarily selected, deterministic initial system state. Then the state evolution sequence $\vn$ generated by either a synchronous or asynchronous implementation of the filtered system operator $T_f$ with $\rho\in(0,1)$ satisfies $\v^n \xrightarrow{(2,2)} \v^\star$ where $\v^\star$ is a fixed-point of $T$.
		\end{lem}

		We now provide commentary on the interpretation of several steps of the proof of Lemma~\ref{lem:passive-everywhere} in order to underscore the geometric reasons why convergence guarantees require a filtered system operator. For clarity, we assume in our discussion and the accompanying illustrations that the fixed-point set $\mathcal{F}_{T}$ is singleton and that the distances from both $\v$ and $T(\v)$ to $\mathcal{F}_{T}$ are equal, i.e.~we assume the ``worst case'' scenario consistent with achieving the upper bound of the inequality:
		\begin{eqnarray}
			\|T(\v) - \v^\star\|^2 = \|T(\v) - T(\v^\star)\|^2 \leq \|\v - \v^\star\|^2.
		\end{eqnarray}
		The key identity used in establishing the inequalities to be iterated for both synchronous and asynchonous arguments is:
		\begin{eqnarray} \label{eq:pre-stewart}
			\left\|\rho T\left(\v^{n-1}\right) + (1-\rho)\v^{n-1} - \v^\star\right\|^2 \,\, = \,\, \rho \left\|T\left(\v^{n-1}\right) - \v^\star \right\|^2 + (1-\rho)\left\|\v^{n-1} - \v^\star\right\|^2 - \rho(1-\rho)\left\|T\left(\v^{n-1}\right) - \v^{n-1}\right\|^2.
		\end{eqnarray}
		This equality can be understood geometrically by application of Stewarts theorem \cite{stewart} to the filtered system update procedure and is illustrated in Fig.~\ref{fig:passive-everywhere-triangle} where $\v^{n-1} = \v$ and $\v^{n} = \rho T(\v) + (1-\rho)\v$. Specifically, consider a triangle with sides of length $a$, $b$, and $c$ and let $d$ be the length of a cevian to side $a$ such that $d$ divides $a$ into two pieces of lengths $n$ and $m$ where $m$ is adjacent to $c$ and $n$ is adjacent to $b$. Then, Stewart's theorem ensures that these distances satisfy $b^2m+c^2n = a(d^2+mn)$ which we rewrite as:
		\begin{eqnarray} \label{eq:stewart}
			d^2 = \frac{m}{a}b^2 + \frac{n}{a}c^2 - mn. 
		\end{eqnarray}
		By proper assignment of the vertices of the triangle (depicted on the right of Fig.~\ref{fig:passive-everywhere-triangle}), \eqref{eq:pre-stewart} follows immediately from \eqref{eq:stewart}. To this point, however, we did not need to utilize any system operator properties of $T$. In Fig.~\ref{fig:passive-everywhere-triangle} this relationship ensure that $b\leq c$. Indeed, while the proof of Lemma~\ref{lem:passive-everywhere} requires the use of a filtered system operator, in practice we repeatedly observe that a direct asynchronous implementation of $T$ is sufficient for convergence (at an approximately linear rate).  Consistent with the discussion following Definition~\ref{def:asynchronous}, this is in part due to the implicit filtering that occurs in accordance with the stochastic dynamics of $D^{(p)}$; we elaborate on this filtering next. 
		\begin{figure}[t!]
			\centering
		  	\centerline{\includegraphics[width=\textwidth]{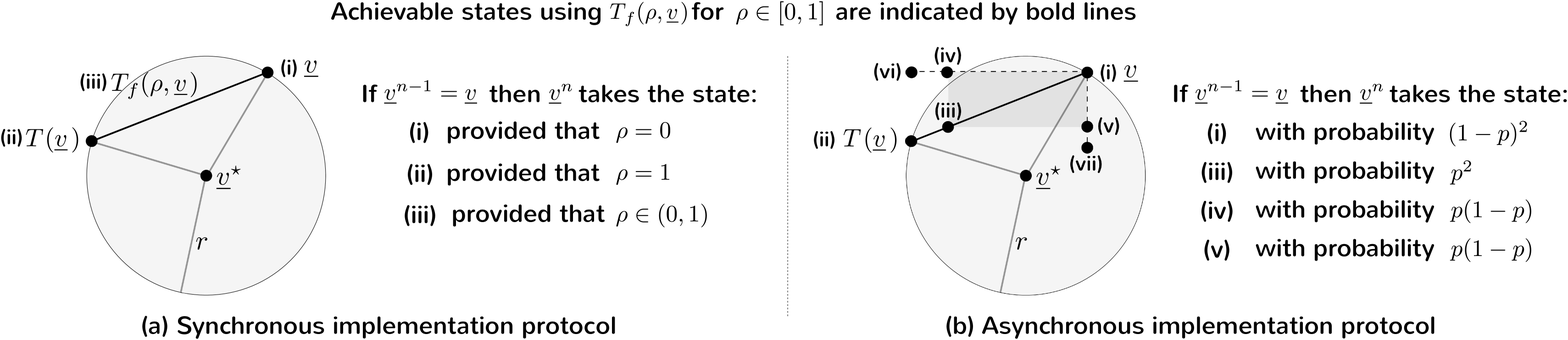}}\caption{An illustration of the geometry underlying the convergence guarantees established in Lemma~\ref{lem:passive-everywhere}. Specifically, setting $\v^{n-1} = \v$, the possible outcomes for $\v^{n}$ are separately depicted for (a) a synchronous implementation and (b) an asynchronous implementation of a filtered system operator where the original system operator $T\colon \R^2 \to \R^2$ is passive everywhere and $\mathcal{F}_{T}$ is assumed to be singleton for clarity.}
			\label{fig:passive-filtering}
		\end{figure}

		A geometric illustration of the update mechanism using a filtered system operator in accordance with Lemma~\ref{lem:passive-everywhere} is provided in Fig.~\ref{fig:passive-filtering} for a passive everywhere system operator $T:\R^2 \to \R^2$. Referring to both implementation protocols (a) and (b) in Fig.~\ref{fig:passive-filtering}, the labels (i)-(iii) indicate: (i) an arbitrary system state $\v$ on the boundary of $B(\v^\star,r)$, (ii) the result of applying $T$ to $\v$, and (iii) the continuum of system states $\rho T(\v) + (1-\rho)\v$ achievable by evaluating $T_{f}$ for $\rho \in [0,1]$. Setting $\v^{n-1} = \v$, we now enumerate the possible states taken by $\v^n$ for each implementation protocol separately. In the synchronous setting two outcomes are generally possible:  (1) the action of $T$ on $\v$ yields $\v$ in which case $\v^n$ is itself a fixed-point of $T$, or (2) the state $\v^n$ is at most equidistant from $\v^\star$ and therefore every state of $T_{f}$ is strictly closer to $\v^\star$ than $\v$ is provided that $\rho \in (0,1)$. Note that sublinear convergence is in part a consequence of the fact that $T(\v)$ may be arbitrarily close to $\v$ on the perimeter of $B(\v^\star,r)$. In the asynchronous setting, the possible states taken by $\v^n$ for $\rho \in [0,1]$ correspond to the states along three line segments: (1) the (solid) line segment with endpoints (i) and (ii), (2) the (dashed) line segment with the endpoints (i) and (vi), and (3) the (dashed) line segment with endpoints (i) and (vii). For a fixed value of $\rho \in (0,1)$, only four outcomes are possible corresponding to labels (i) and (iii)-(v) where the probability of each is listed in the figure. In the general case, the $2^k$ corners of the $k$-hypertope (hyperrectangle) defined with opposite verticies given by $\v$ and $T_f(\rho,\v)$ enumerate the possible outcomes for $\v^n$, as is indicated by the shaded rectangle in the figure. Note that in the asynchronous setting one outcome for $\v^n$ is strictly further than $\v^{n-1}$ from $\v^\star$ despite the convergence guarantee in Lemma~\ref{lem:passive-everywhere}. This is consistent with the previously addressed comparison between almost sure and mean square convergence and further underscores the fact that almost sure convergence is not precluded under the assumptions made but would require further conditions on $T$.  For both (a) and (b), we see that the optimal selection of $\rho$ in the sense of a provable convergence rate is $\frac{1}{2}$ since it minimizes the bound in \eqref{eq:passive-convergence-rate} for any $p\in(0,1]$. Graphically, this corresponds to achieving the state along the line segment with endpoints (i) and (ii) which is closest to $\v^\star$. Analytically, this is equivalent to selecting the state which is the orthogonal projection of $\v^\star$ onto the affine subspace defined by $\rho T(\v) + (1-\rho)\v$ for all $\rho \in \R$.

	\subsection{Convergence and stability results related to general system operators}  \label{sec:general}
		We now consider the class of general system operators and in particular those which are associated with $\alpha$-expansivity in some sense. Building upon the discussion in Section~\ref{sec:passive}, the set of fixed-points $\mathcal{F}_{T}$ associated with an $\alpha$-expansive everywhere system operator $T$ can be arbitrary and thus need not be convex when non-empty. We proceed in this section with the minimal assumption that $\mathcal{F}_{T}$ is non-empty. 
		
		Sufficient conditions for the convergence of a synchronous or asynchronous implementation of an $\alpha$-expansive everywhere system operator to a fixed-point cannot generally be established. However, continuing our discussion in a Hilbert space, it is possible to make judicious use of a filtered system operator in order to ensure convergence so long as the system operator has a mixing property resembling that of topological mixing in chaos theory \cite{chaos-theory}. We state this more formally in the next lemma.

		\begin{lem}	\label{lem:expansive-everywhere}
			Let $T$ denote a system operator which is $\alpha$-expansive everywhere with a non-empty convex set of fixed-points $\mathcal{F}_{T}$ and let $\v^0$ be an arbitrarily selected, deterministic initial system state. If $T$ additionally satisfies the mixing property:
			\begin{eqnarray}
				\sup_{\v\not\in\mathcal{F}_{T}}\frac{\left\langle T(\v)-\v^\star,\v-\v^\star \right\rangle}{\left\|T(\v)-\v^\star\right\|\left\|\v-\v^\star\right\|} \leq \gamma \label{eq-lem:expansive-everywhere-mixing-condition}
			\end{eqnarray} 
			for all $\v^\star \in \mathcal{F}_{T}$ and some $\gamma\in[-1,1)$ such that $\alpha \gamma < 1$, then the state evolution sequence $\vn$ associated with either a synchronous or asynchronous implementation of the filtered system operator $T_f$ with:
			\begin{eqnarray}
				\rho \in \left( 0, \frac{2(1-\alpha\gamma)}{1 + \alpha^2 - 2\alpha\gamma} \right) \label{eq-lem:expansive-everywhere-filter-parameter}
			\end{eqnarray}
			satisfies $\v^n \xrightarrow{(1,2)} \v^\star$ where $\v^\star$ is a fixed-point of $T$.
		\end{lem}

		\begin{figure}[t!]
			\centering
		  	\centerline{\includegraphics[width=\textwidth]{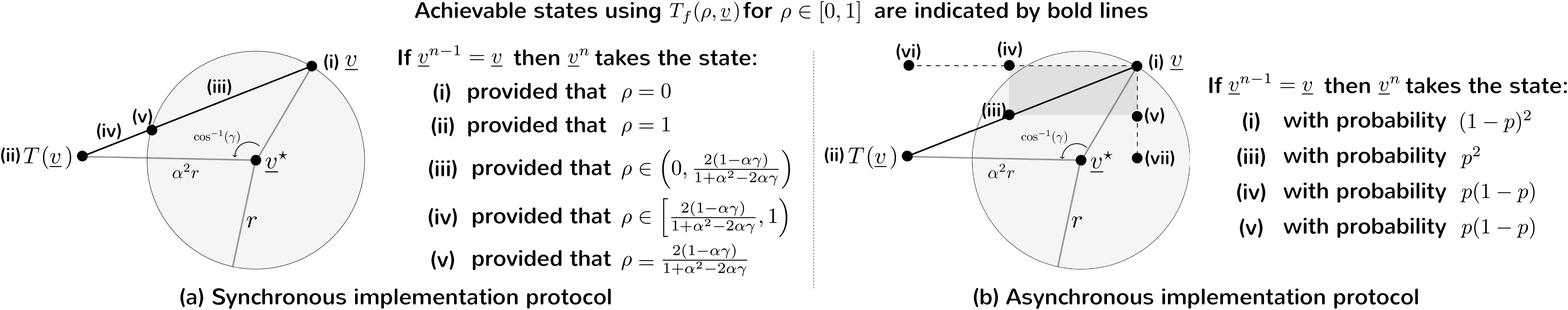}}\caption{An illustration of the geometry underlying the convergence guarantees established in Lemma~\ref{lem:expansive-everywhere}. Specifically, setting $\v^{n-1} = \v$, the possible outcomes for $\v^{n}$ are separately depicted for (a) a synchronous implementation and (b) an asynchronous implementation of a filtered system operator where the original system operator $T\colon \R^2 \to \R^2$ is $\alpha$-expansive everywhere and $\mathcal{F}_{T}$ is assumed to be singleton for clarity. The role of $\gamma$ in \eqref{eq-lem:expansive-everywhere-mixing-condition} is also indicated.}
			\label{fig:expansive-filtering}
		\end{figure}

		A geometric illustration of the convergence properties summarized by Lemma~\ref{lem:expansive-everywhere} is provided in Fig.~\ref{fig:expansive-filtering} for an $\alpha$-expansive everywhere system operator $T:\R^2 \to \R^2$ that additionally satisfies \eqref{eq-lem:expansive-everywhere-mixing-condition}. For clarity, we assume in the illustration that the fixed-point set $\mathcal{F}_{T}$ is singleton and that the distance from $T(\v)$ to $\mathcal{F}_{T}$ is maximized with respect to the distance from $\v$ to $\mathcal{F}_{T}$, i.e.~we assume the ``worst case'' scenario consistent with achieving the upper bound of the inequality 
		\begin{eqnarray}
			\|T(\v) - \v^\star\|^2 = \|T(\v) - T(\v^\star)\|^2 \leq \alpha^2\|\v - \v^\star\|^2.
		\end{eqnarray}
		Referring to both implementation protocols (a) and (b) in Fig.~\ref{fig:expansive-filtering}, the label (i) indicates an arbitrary system state $\v$ on the boundary of $B(\v^\star,r)$ while the label (ii) indicates the result of applying $T$ to $\v$. Moreover, the role of the mixing parameter $\gamma$ in \eqref{eq-lem:expansive-everywhere-mixing-condition} is indicated as the cosine of the angle between the vectors pointing from $\v^\star$ to $\v$ and from $\v^\star$ to $T(\v)$. In the synchronous setting the labels (iii)-(v) indicate: (iii) the continuum of system states strictly closer to $\v^\star$ than $\v$ achievable by evaluating $T_{f}$ for $\rho$ satisfying \eqref{eq-lem:expansive-everywhere-filter-parameter}, (iv) the system states along the same trajectory equal to or further from $\v^\star$ than $\v$, and (v) the system state equidistant with $\v$ from $\v^\star$ corresponding to evaluating $T_{f}$ for $\rho = \frac{2(1-\alpha\gamma)}{1 + \alpha^2 - 2\alpha\gamma}$. By proper selection of the filter parameter $\rho$ it is evident from the figure that we are able to ensure the filtered system operator acts as a dissipative everywhere system operator since $\gamma < 1$ guarantees that $\v$ and $T(\v)$ cannot be arbitrarily close together; this conclusion is stronger than the conclusion made in the discussion surrounding Fig.~\ref{fig:passive-filtering}. In the asynchronous setting, the possible states taken by $\v^n$ when $\v^{n-1}=\v$ as a function of $\rho \in [0,1]$ correspond to the states along three line segments: (1) the (solid) line segment with endpoints (i) and (ii), (2) the (dashed) line segment with the endpoints (i) and (vi), and (3) the (dashed) line segment with endpoints (i) and (vii). For a fixed value of $\rho$, only four outcomes are possible corresponding to labels (i) and (iii)-(v) where the probability of each is listed in the figure. Analagous to the discussion in Section~\ref{sec:passive}, the $2^k$ corners of the $k$-hypertope (hyperrectangle) defined with opposite verticies given by $\v$ and $T_f(\rho,\v)$ enumerate the possible outcomes, and is again depicted by the shaded rectangle in the figure for $k=2$. Note that in the asynchronous setting it is possible for $\v^n$ to be strictly further than $\v^{n-1}$ from $\v^\star$ both due to the stochastic operator $D^{(p)}$ and the filter coefficient $\rho$. The optimal selection of $\rho$ in the sense of a provable convergence rate, obtained in the proof of Lemma~\ref{lem:expansive-everywhere}, is given by:
		\begin{eqnarray} \label{eq:rate-expansive}
			\rho = \frac{1 - \alpha\gamma}{1+\alpha^2 - 2\alpha\gamma}.
		\end{eqnarray}
		Graphically, this corresponds to achieving the state along the line segment with endpoints (i) and (iii) which is closest to $\v^\star$.  Analytically, this is equivalent to selecting the state which is the orthogonal projection of $\v^\star$ onto the affine subspace defined by $\rho T(\v) + (1-\rho)\v$ for all $\rho \in \R$.

		We now briefly comment on the relationship between the class of system operators which satisfy Lemma~\ref{lem:expansive-everywhere} and those which are chaotic. Recall that a given system operator $T$ is chaotic if the state evolution sequence produced using a synchronous implementation protocol satisfies three conditions: (1) hypersensitivity to initial conditions, (2) topological mixing, and (3) dense periodic orbits. While $\alpha$-expansive everywhere system operators generally satisfy (1), this is not sufficient to be chaotic. For example, consider the system operator $T(v) = -1.1v$ which satisfies both condition (1) and the qualifiers of Lemma~\ref{lem:expansive-everywhere} but is not chaotic. Verifying hypersensitivity to initial conditions for this system operator follows from defining the state evolution sequences $\vn$ and $\un$ respectively initialized with $\v^0$ and $\u^0 = \v^0+\delta$ for some $\delta \neq 0$ and observing that:  
		\begin{eqnarray}
			\v^n = (-1.1)^n\v^0 & \hspace{1em} \text{ and } \hspace{1em} & \u^n = (-1.1)^n\left(\v^0 + \delta\right)
		\end{eqnarray}
		hence $\left\| \v^n - \u^n \right\| = (1.1)^n|\delta|$ can be made arbitrarily large by proper selection of $n\in\N$. While Lemma~\ref{lem:expansive-everywhere} requires mixing in the sense of \eqref{eq-lem:expansive-everywhere-mixing-condition}, topological mixing for chaotic operators is stronger since $T$ is required to be hypercyclic about some point, i.e.~$\vn$ must be dense in $\Rk$. Finally, Lemma~\ref{lem:expansive-everywhere} makes no requirement similar to condition (3).

		\begin{figure}[t]
			\begin{minipage}{0.45\textwidth}
				\centerline{\includegraphics[width=\textwidth]{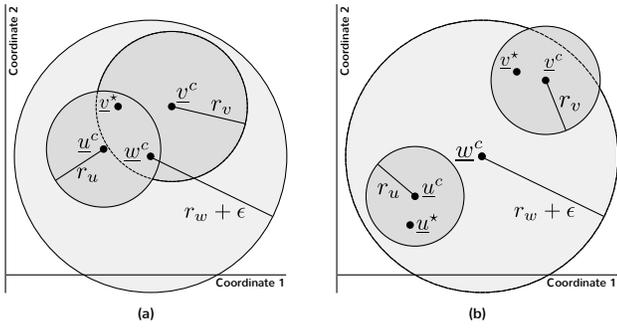}}
			\end{minipage} \hfill
			\begin{minipage}{0.55\textwidth}
				\normalsize 
				Let $T:\R^{2} \to \R^{2}$ denote a system operator satisfying the following: 
				\begin{itemize}
					\item[(i)] Property: $T$ is $\alpha_w$-dissipative about $\w^c$ where $\w^c \not \in \mathcal{F}_{T}$. 
					\item[] Reference: Lemma~\ref{lem:dissipative-about-some-point}
					\vspace{.1in}
					\item[(ii)]	Property: $T$ is $\alpha_v$-dissipative about all $\v\in B(\v^c,r_v)$ for some $r_v$ satisfying $\|T(\v^c) - \v^c\| \leq r_v(1-\alpha_v)$. 
					\item[] Reference: Lemma~\ref{lem:dissipative-everywhere-ball} 
					\vspace{.1in}
					\item[(iii)] Property: $T$ is $\alpha_u$-dissipative about all $\u\in B(\u^c,r_u)$ for some $r_u$ satisfying $\|T(\u^c) - \u^c\| \leq r_u(1-\alpha_u)$.
					\item[] Reference: Lemma~\ref{lem:dissipative-everywhere-ball}
				\end{itemize}
			\end{minipage}	
			\caption{\label{fig:general-system-operator-example} An example of the type of aggregate convergence guarantees that can be made by combining the Lemmas presented in Section~\ref{sec:dissipative}. Several Euclidean balls in the codomain of a general system operator $T \colon \R^2 \to \R^2$ are illustrated on the left corresponding to the properties listed on the right.}
		\end{figure}
		
		Finally, we illustrate the utility of the lemmas presented hereto in establishing convergence and fixed-point properties of an arbitrary system operator $T$ using the example in Fig.~\ref{fig:general-system-operator-example}. Referring to this figure, the system operator $T:\R^2 \to \R^2$ is assumed to satisfy the properties listed on the right whose consequences are summarized by the companion lemmas, also listed. On the left are two possible scenarios consistent with these assumptions. From property (i) we conclude that the state evolution sequence $\vn$ will eventually (in finite time) be contained within the Euclidean ball $B(\w^c, r_w+\epsilon)$ for any $\epsilon > 0$ regardless of the initial system state $\v^0$. From properties (ii) and (iii) we respectively conclude that there exists a unique fixed-point $\v^\star$ of $T$ in $B(\v^{c},r_v)$ and a unique fixed-point $\u^\star$ in $B(\u^{c},r_u)$. Additionally, we conclude that if the state evolution sequence ever enters $B(\v^c,r_v)$ then it will converge linearly to $\v^\star$ and likewise that if the state evolution sequence ever enters $B(\u^c,r_u)$ it will converge linearly to $\u^\star$. These conclusions are true for both scenarios. Scenario (a) depicts the case where $B(\u^c, r_u) \cap B(\v^c, r_v) \neq \emptyset$ and therefore $\v^\star = \u^\star$. We comment that if the state evolution sequence enters the region $B(\u^c, r_u) \cap B(\v^c, r_v)$ then it will never leave it. If $T$ additionally was $\alpha$-disispative about $\v^\star = \u^\star$ (Lemma~\ref{lem:dissipative-about-fixed-point}) then we would be able to conclude convergence to this fixed-point starting from any initial state. Scenario (b) depicts the case where $B(\u^c, r_u) \cap B(\v^c, r_v) = \emptyset$ hence the system operator $T$ will have two distinct fixed-points. So, if state evolution sequence $\vn$ enters either $B(\u^c, r_u)$ or $B(\v^c, r_v)$ it will converge linearly to the unique fixed-point inside that ball.

\section{Numerical experiments} \label{sec:examples}	
	In this section we investigate the numerical convergence properties of synchronous and asynchronous implementations of system operators associated with various conservative constraint satisfaction problems. In order to do this, we perform a series of experiments outlined as follows: 
	\begin{itemize}
		\item[1] for a given dimensionality $k$, realize the system operator $T:\Rk\to\Rk$ associated with the CCSP at hand; 
		\item[2] identify a fixed-point $\v^\star \in \mathcal{F}_{T}$ of the system operator;
		\item[3] for an appropriate implementation of $T$, generate the state evolution sequence $\vn$ for various values of $p$ and track the distance of $\v^n$ to $\v^\star$;
		\item[4] repeat 3 for many trials and average the results.
	\end{itemize}
	The result of this procedure is depicted in each subsection below as a function of ``equivalent (normalized) iterations'' by which we mean the iteration count corresponding to the total computation performed. This quantity is approximately given by the index $n$ times the probability $p$ associated with the stochastic matrix $D^{(p)}$.

		\begin{figure}[t!]
			\centering
		  	\centerline{\includegraphics[width=\textwidth]{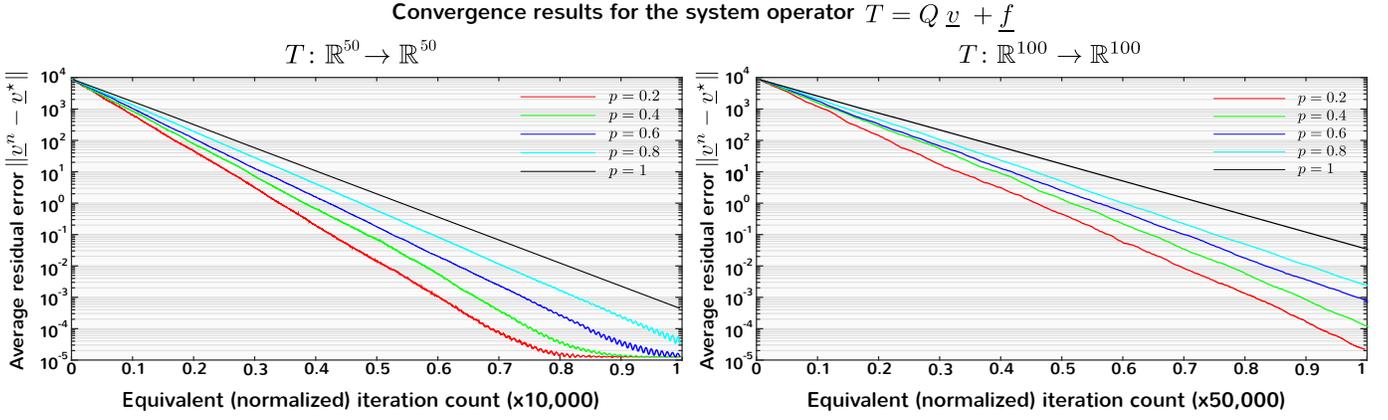}}\caption{Numerical convergence behavior of a filtered system implementation of the passive everywhere source system operator \eqref{eq:example-passive-everywhere} with filter coefficient $\rho = 0.5$. The performance depicted is averaged over $1000$ trials for various values of $p$ and dimensionalities $k=25$ (left) and $k=100$ (right). The initial system state $\v^0$ is randomly selected on the boundary of $B(\underline{0},1)$ for each trial. Lemma~\ref{lem:passive-everywhere} summarizes the theoretical guarantees for this case.}
			\label{fig:example-passive-everywhere}
		\end{figure}
	\subsection{Convergence of a passive everywhere source system operator}
		We now consider the class of CCSPs embodied by a source system operator as described by Definition~\ref{def:source}. Specifically, let $T \colon \Rk \to \Rk$ denote a passive everywhere source system operator of the form: 
		\begin{eqnarray}	 \label{eq:example-passive-everywhere}
			T(\v) = Q\v + \f
		\end{eqnarray}
		where $Q$ is a random orthogonal matrix obtained by first generating a $k \times k$ matrix from a Gaussian ensemble and projecting it to the nearest orthogonal matrix in the Frobenious norm sense and $\f$ is a Gaussian random vector. Observe that if $\f = \underline{0}$ then \eqref{eq:example-passive-everywhere} reduces to a neutral system operator (Definition~\ref{def:neutral}). For this experiment, we generate the state evolution sequences using a filtered system implementation with filter parameter $\rho = 0.5$, i.e.~using the signal-flow structure in row (3) of column (b) in Fig.~\ref{fig:composition-definitions}. The results of this experiment for dimensionalities $k=50$ and $k=100$ and asynchronous delay probabilities $p=0.2, 0.4, 0.6, 0.8, 1$ are depicted in Fig.~\ref{fig:example-passive-everywhere}. Observe that on average the state evolution sequences initially converge at a linear rate toward the fixed-point and slow to a sublinear rate once they reach a sufficiently nearby region, consistent with the worst-case performance established by Lemma~\ref{lem:passive-everywhere}. As expected, the synchronous implementation ($p=1$) exhibits monotonic convergence. Finally, we comment that using the implementation techniques presented in this paper with \eqref{eq:example-passive-everywhere} yields an asynchronous algorithm for solving linear systems of the form $(I-Q)\v =\f$ without explicitly decomposing or inverting $(I-Q)$ or $Q$. Further, by proper selection of $\f$ a subset of the eigenvectors and eigenvalues of $Q$ may be identified.

		\begin{figure}[t!]
			\centering
		  	\centerline{\includegraphics[width=\textwidth]{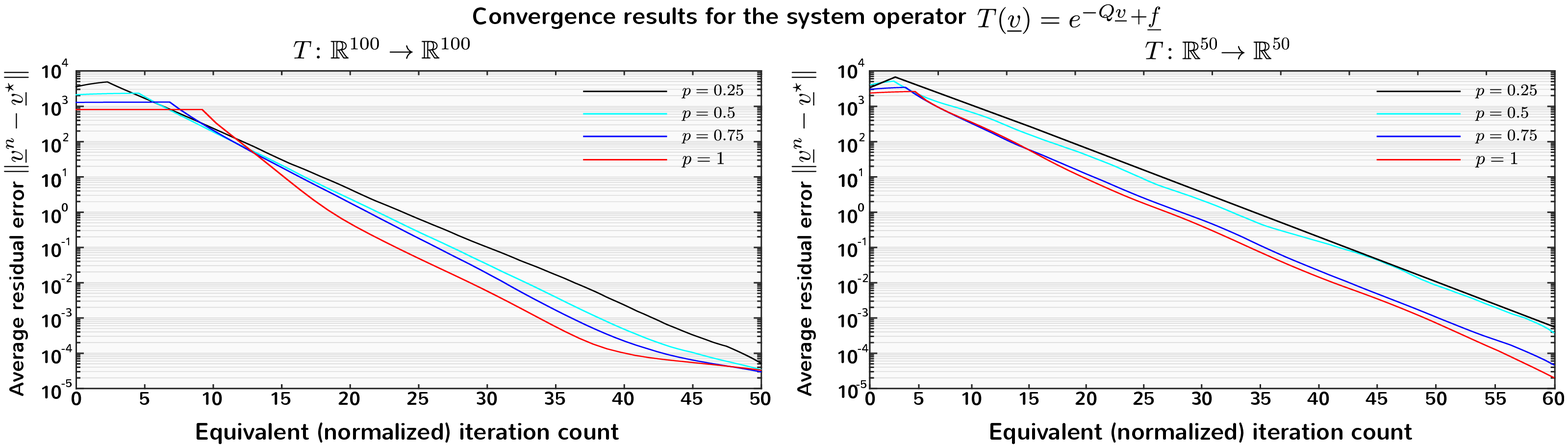}}\caption{Numerical convergence behavior of a direct implementation of the system operator \eqref{eq:example-non-alpha-conic} which does not satisfy \eqref{eq:alpha-conic} for any finite value of $\alpha$. The performance depicted is averaged over $1000$ trials for various values of $p$ and dimensionalities $k=100$ (left) and $k=50$ (right). The initial system state $\v^0$ is randomly selected on the boundary of $B(\underline{0},1)$ for each trial. Lemma~\ref{lem:dissipative-everywhere-composition} summarizes the theoretical guarantees for this case.}
			\label{fig:example-non-alpha-conic}
		\end{figure}
	\subsection{Convergence of a non $\alpha$-conic system operator}
		In the discussion following Lemma~\ref{lem:dissipative-everywhere-composition} we examined the scalar system operator $T(v) = e^{-v}$ and concluded that while $T$ is itself not $\alpha$-conic for any finite value of $\alpha$, the system operator $T\circ T(v)$ is with coefficient $\alpha = \frac{1}{e}$. In this subsection we consider a generalization of this example. Specifically, let $T \colon \Rk \to \Rk$ denote a general system operator of the form:
		\begin{eqnarray}	 \label{eq:example-non-alpha-conic}
			T(\v) = e^{-Q\v} + \f
		\end{eqnarray}
		where $Q$ is a random orthogonal matrix with eigenvalues bounded away from $-1$, $\f$ is a Gaussian random vector, and the exponential is coordinatewise. $Q$ is obtained as follows:  we generate a candidate matrix by first drawing from a Gaussian ensemble and then projecting to the nearest orthogonal matrix in the Frobenious norm sense and continue regenerating candidates until one satisfies the stated eigenvalue property. For this experiment, we generate the state evolution sequences by implementing $T$ directly, i.e.~using the signal-flow structure in row (3) of column (a) in Fig.~\ref{fig:composition-definitions}. The results of this experiment for dimensionalities $k=50$ and $k=100$ and asynchronous delay probabilities $p=0.25, 0.5, 0.75, 1$ are depicted in Fig.~\ref{fig:example-non-alpha-conic}. Observe that after finite-time transient effects terminate, the average performance indicates that the state evolution sequences converge linearly to the fixed-point with a rate that does not depend upon $p$, consistent with the behavior expected from application of Lemma~\ref{lem:dissipative-everywhere-composition} to \eqref{eq:example-non-alpha-conic}. We comment that using the presented implementation techniques with \eqref{eq:example-non-alpha-conic} yields an asynchronous algorithm for solving transcendental exponential equations that have no analytic solution, i.e.~which cannot be expressed in closed form as a polynomial equation.

		\begin{figure}[t!]
			\centering
		  	\centerline{\includegraphics[width=\textwidth]{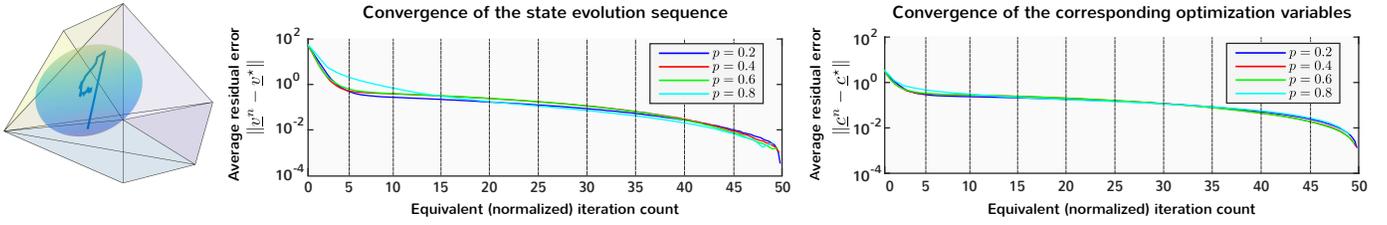}}\caption{Numerical convergence behavior of a direct implementation of the passive everywhere system operator embodying the Chebychev center problem \eqref{eq:example-chebychev-center} as discussed in \cite{LinProg}. The sphere obtained in $\R^3$ on the left illustrates the convergence of the center of the sphere as the associated state evolution sequence converges. The performance depicted in the convergence plots in the middle (state evolution variables $\v^n$) and on the right (optimization variables $\c^n$) are averaged over $500$ trials for various values of $p$ and dimensionality $k=100$. The initial system state $\v^0 = \underline{0}$ for each trial.}
			\label{fig:example-chebychev-center}
		\end{figure}		
	\subsection{Convergence of the Chebychev center problem}
		The Chebychev center problem is a linear programming problem used to identify the largest Euclidean ball $B(\c,r)$ which is contained within a closed convex polytope. In particular, if the polytope is characterized in half-space notation by $\{\v\in\R^k \colon A\v \leq \b \}$ for a given matrix $A\in\R^{m\times k}$ and vector $b\in\R^m$, then the Chebychev center is identified by solving: 
		\begin{eqnarray}\label{eq:example-chebychev-center}
			\begin{array}{rl}
			\displaystyle \min_{\c, r} & -r \\
			\text{s.t.} & \langle \a^i,\c\rangle + \| \a^i\| r  \leq  \b_i, \hspace{1em} i = 1, \dots , m\\
			& r \geq 0
			\end{array}
		\end{eqnarray}
		where $\a^i$ is the $i$th column of $A$. The general methodology behind posing optimization problems as conservative constraint satisfaction problems is discussed in \cite{BLGS1} and the specific formulation of the system operator $T$ for this problem is provided in \cite{LinProg}.  Without loss of generality, the system operator for this problem as well as general linear programs can be shown to be passive everywhere. We note here for context that the variables encapsulated by the state evolution sequence vector $\v^n$ are to within a change of basis the optimization variables hence the implicit state evolution sequences $\{\c^n\}_{n=0}^{\infty}$ and $\{r^n\}_{n=0}^{\infty}$ are defined for a given state evolution sequence $\vn$ associated with implementing $T$. For this experiment, we generate the state evolution sequences by implementing $T$ directly, i.e.~using the signal-flow structure in row (3) of column (a) in Fig.~\ref{fig:composition-definitions}. Our choice to implement this system operator directly is specifically to illustrate the previously addressed observation that asynchronous implementations of passive everywhere system operators  in practice typically converge due to implicit filtering although such a guarantee is not established in this paper. The results of this experiment for a randomly generated convex polytope with $m=200$ faces in $k=100$ dimensions and asynchronous delay probabilities $p=0.2, 0.4, 0.6, 0.8$ are depicted in Fig.~\ref{fig:example-chebychev-center}. Observe that the average performance indicates sublinear convergence, consistent with the convergence expected had we implemented the system operator using a filtered implementation.

		\begin{figure}[t!]
			\centering
		  	\centerline{\includegraphics[width=\textwidth]{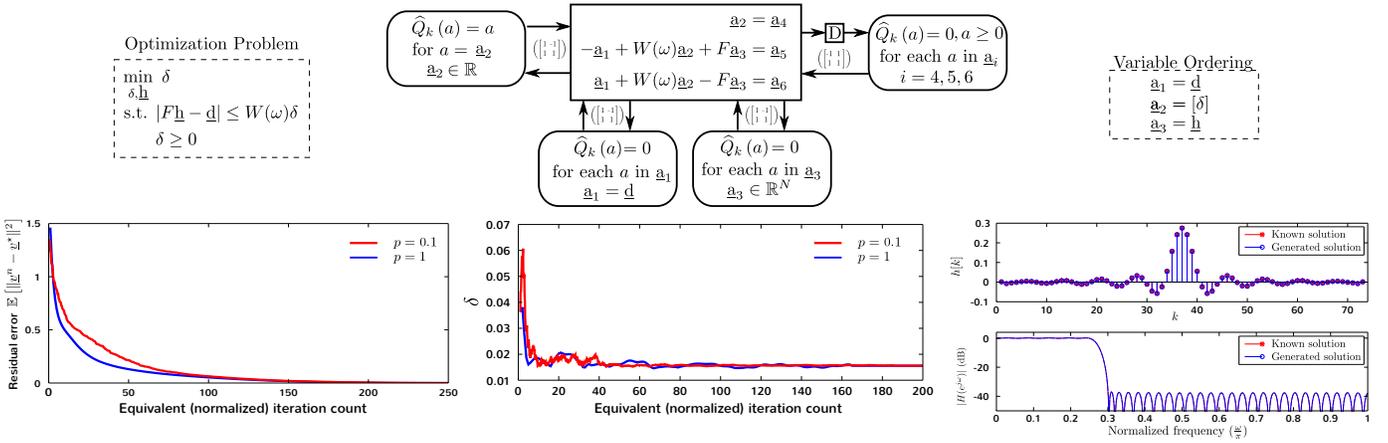}}\caption{The conservative signal-flow structure and numerical convergence behavior of a direct implementation of the passive everywhere system operator embodying the minimax optimal filter design problem \eqref{eq:example-filter-design} as discussed in \cite{BLGS2}. The performance depicted is taken over a single trial for $p=0.1$ and $p=1$ and dimensionality $k=37$ with initial system state $\v^0 = \underline{0}$.}
			\label{fig:example-filter-design-1}
		\end{figure}
	\subsection{Convergence of the minimax optimal filter design problem}
\begin{figure}[t!]
			\centering
		  	\centerline{\includegraphics[width=\textwidth]{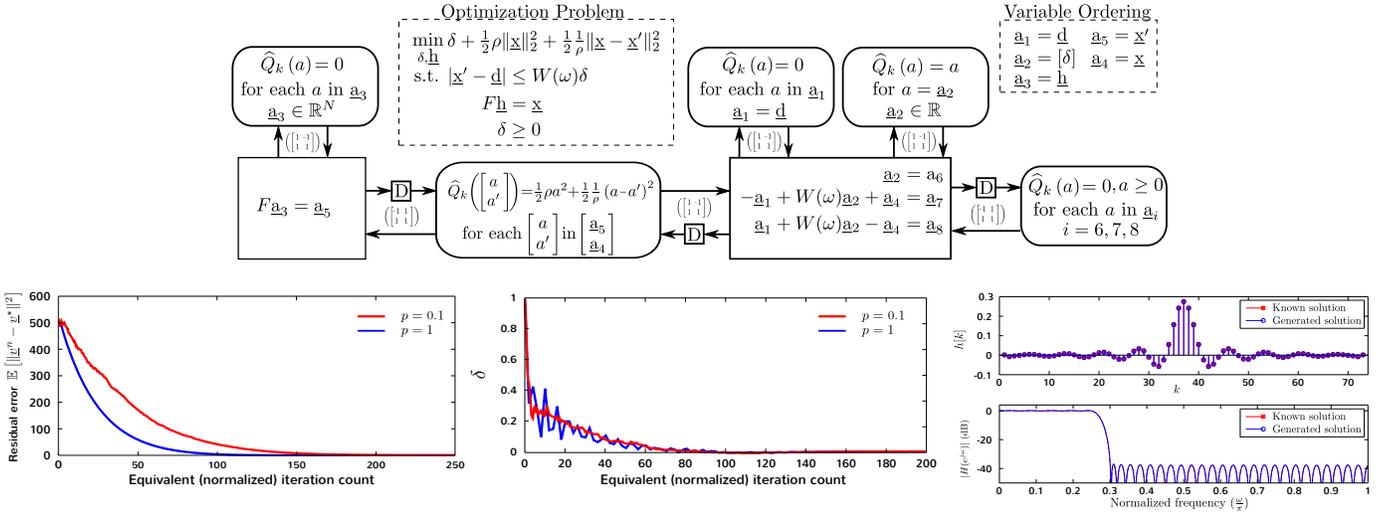}}\caption{The conservative signal-flow structure and numerical convergence behavior of a direct implementation of the passive everywhere system operator embodying the formulation of the minimax optimal filter design problem as a convex quadratic program as discussed in \cite{BLGS2}. The performance depicted is taken over a single trial for $p=0.1$ and $p=1$ and dimensionality $k=37$ with initial system state $\v^0 = \underline{0}$.}
			\label{fig:example-filter-design-2}
		\end{figure}

		The design of a minimax optimal filter is a classical problem in signal processing and is often solved using the Parks-McClellan algorithm, a modification of the more general Remez exchange algorithm. Alternatively, the design of a symmetric, finite impulse response filter with support $[0,2q]$ characterized by $\h\in\R^{q+1}$ given a target frequency response $\d\in\R^{m}$ according to the minimax criterion is also obtained by solving the linear program: 
		\begin{eqnarray} \label{eq:example-filter-design}
			\begin{array}{rl}
			\displaystyle \min_{\delta, \h} & \delta \\
			\text{s.t.} & |F\h - \d| \leq W(\omega) \delta \\
			& \delta \geq 0
			\end{array}
		\end{eqnarray}
		where $F\in\R^{m\times (q+1)}$ maps the filter coefficients $\h$ from the sample domain to the frequency domain, $W(\omega)\in\R^{m}$ is a positive scaling vector and $\delta\in\R$ is the achieved deviation. We perform two experiments for the design of this filter. The conservative signal-flow structures associated with the passive everywhere system operators for each experiment are respectively provided in Fig.~\ref{fig:example-filter-design-1} and Fig.~\ref{fig:example-filter-design-2} and have been adapted from \cite{BLGS2}.  For the first experiment, the state evolution sequences are generated by directly implementing the system operator $T$ embodying \eqref{eq:example-filter-design}, i.e.~using the signal-flow structure in row (3) of column (a) in Fig.~\ref{fig:composition-definitions}.  For the second experiment, we again generate the state evolution sequences by directly implementing the system operator $T$ embodying the convex quadratic formulation of the design problem listed in Fig.~\ref{fig:example-filter-design-2}. The results for each of these experiments with $m=1000$ and $k=37$ are depicted in Figs.~\ref{fig:example-filter-design-1} and \ref{fig:example-filter-design-2}, respectively, for asynchronous delay probabilities $p=0.1$ and $p=1$. Similar to the previous subsection, our choice to implement the system operators for both experiments directly is again to illustrate the convergence of a passive everywhere system operator in practice. Observe that the average rate of convergence seen in both experiments is linear and is better than the rate experienced in the previous subsection for the same class of passive everywhere system operators.

\cleardoublepage\FloatBarrier
\appendix

	\begin{pf}[Lemma~\ref{lem:translate-to-zero}] 
		We prove this result by first providing the explicit form of one such system operator $T_{(2)}$ and then verifying that the properties claimed do in fact hold. Define the system operator $T_{(2)} \colon \Rk \to \Rk$  as: 
		\begin{eqnarray}
			T_{(2)}(\v) \triangleq T_{(1)}(\v + \v^1)- \v^1. \label{eq-lem:translate-to-zero-10}
		\end{eqnarray}
		Substituting $T_{(2)}$ into the $\alpha$-conic about $\underline{0}$ condition from Definition~\ref{def:alpha-conic} yields:
		\begin{eqnarray}
			\sup_{\u \neq \underline{0}} \frac{\| T_{(2)}(\underline{0} + \u) - T_{(2)}(\underline{0})\|}{\| \u \|}  &=&  \sup_{\u \neq \underline{0}} \frac{\| T_{(1)}(\u + \v^1) - T_{(1)}(\v^1)\|}{\| \u \|} \\ 
			&\leq&  \alpha
		\end{eqnarray} 
		where the inequality 
		is due to the assumption that $T_{(1)}$ is $\alpha$-conic about $\v^1$. Therefore we conclude that $T_{(2)}$ is indeed $\alpha$-conic about $\underline{0}$ with the same parameter $\alpha$. The map relating the fixed-points in $\mathcal{F}_{T_{(1)}}$ and $\mathcal{F}_{T_{(2)}}$ corresponds to an affine translation by $\v^1$ and follows from straightforward algebra. For the particular choice of system operator $T_{(2)}$ in \eqref{eq-lem:translate-to-zero-10} the fixed-points sets $\mathcal{F}_{T_{(1)}}$ and $\mathcal{F}_{T_{(2)}}$ satisfy $\mathcal{F}_{T_{(2)}} = \mathcal{F}_{T_{(1)}} - \v^1$, i.e.~for each fixed-point $\u^\star \in \mathcal{F}_{2}$ there exists a fixed-point $\v^\star \in \mathcal{F}_{1}$ which satisfies the relationship $\u^\star  = \v^\star - \v^1$. 
	\end{pf}

	\begin{pf}[Lemma~\ref{lem:conic-composition}]
		Before proving this lemma we state two supporting facts. First, by assumption we have that the range of each system operator $T_{(i)}$ is a subset of the domain of $T_{(i+1)}$ for $i = 1, 2, \dots, m-1$ thus validating the necessary domain-codomain intersections of the pairwise composition operators. Second, let $T_{(a)}$ and $T_{(b)}$ be two system operators each satisfying \eqref{eq:alpha-conic-everywhere-equivalent} with respective parameters $\alpha_{a}$ and $\alpha_{b}$. Observe:
		\begin{eqnarray}
			\left\| T_{(b)}\left(  T_{(a)}\left( \v' \right)\right) - T_{(b)}\left( T_{(a)}\left( \u' \right)\right) \right\| & \leq & \alpha_{b}\left\|T_{(a)}\left( \v' \right) - T_{(a)}\left( \u'\right) \right\| \\
			& \leq & \alpha_{a}\alpha_{b}\left\| \v'  -  \u' \right\| \label{eq-lem:conic-composition-10}
		\end{eqnarray}
 		from which we conclude that $T_{(b)}\circ T_{(a)}$ is $\alpha$-conic everywhere with parameter $\alpha_{a}\alpha_{b}$. We now prove the lemma by induction:
 		\begin{itemize}
 			\item[] Base case: The map $T_{(2)} \circ T_{(1)}$ is $\alpha$-conic with parameter $\alpha_1\alpha_2$. This follows immedietely from \eqref{eq-lem:conic-composition-10} by assigning $a=1$ and $b=2$.

 			\item[] Induction step: Suppose that the system operator $T_{(m-1)}\circ\cdots\circ T_{(1)}$ is $\alpha$-conic with parameter $\prod_{i=1}^{m-1}\alpha_{i}$. It then also follows immediately from \eqref{eq-lem:conic-composition-10} that the system operator $T_{(m)} \circ \cdots \circ T_{(1)}$ is $\alpha$ conic with parameter $\prod_{i=1}^{m}\alpha_{i}$ by assigning $T_{(a)} = T_{(m-1)}\circ\cdots\circ T_{(1)}$, $T_{(b)} = T_{(m)}$, $\alpha_{a} = \prod_{i=1}^{m-1}\alpha_{i}$ and $\alpha_{b}= \alpha_m$.
 		\end{itemize}
 		This concludes the proof that the $\alpha$-conic everywhere property is invariant with respect to composition.
	\end{pf}
		
	\begin{pf}[Lemma~\ref{lem:dissipative-about-fixed-point}]
		We begin by proving that the state evolution sequence $\vn$ genereated by a synchronous implementation of $T$ converges to the assumed fixed-point $\v^\star\in\mathcal{F}_{T}$. To this end, define a second state evolution sequence $\un$ in $\Rk$ as the translated by $\v^\star$ system states, i.e. $\u^n = \v^n - \v^\star$ for $n = 0,1,2,\dots$. We then obtain the inequality: 
		\begin{eqnarray}
			\left\|\v^{n} - \v^\star \right\| & =    &   \left\| T(\v^{n-1}) - T(v^\star)\right\|                 \\
											  & =    &   \left\| T(\u^{n-1}+\v^\star) - T(\v^\star)\right\|       \\
											  & \leq &   \alpha \left\| \u^{n-1} \right\|                         \\
											  & =    &   \alpha \left\| \v^{n-1} - \v^\star \right\| \label{eq-lem:dissipative-about-fixed-point-10}
		\end{eqnarray}
		where the first equality is due to Definition~\ref{def:synchronous} and the fact $\v^\star \in \mathcal{F}_{T}$, the second equality is due to the translated system states, the third inequailty is due to the $\alpha$-dissipative about $\v^\star$ assumption, and the final equality is again due to the translated system states. Iterating the inequality in \eqref{eq-lem:dissipative-about-fixed-point-10} $n$ times and taking a limit yields:
		\begin{eqnarray}
			\lim_{n\to \infty} \left\| \v^{n} - \v^\star \right\| &\leq& \left\| \v^0 - \v^\star \right\| \lim_{n\to \infty}\alpha^n \\
			& = & 0.
		\end{eqnarray}
		Since, by assumption, we have that $T(\v^\star)= \v^\star$ and that $T$ is continuous at $\v^\star$ we have proven that $\v^n$ converges to $\v^\star$ for $\ell=1$. We now prove that the state evolution sequence $\vn$ generated by an asynchronous implementation of $T$ also converges to $\v^\star$. Utilizing the translated system states $\u^n = \v^n - \v^\star$ again, we obtain the inequality:
		\begin{eqnarray}
			\E\left[ \left\| \v^n - \v^\star \right\|^2 \right] & = & p \E\left[ \left\|T\left(\v^{n-1}\right) - T\left(\v^\star\right)\right\|^2\right] + (1-p)\E\left[\left\|\v^{n-1} - \v^\star \right\|^2\right] \\ 
			& = & p \E\left[ \left\|T\left(\u^{n-1}+\v^\star\right) - T\left(\v^\star\right)\right\|^2\right] + (1-p)\E\left[\left\|\v^{n-1} - \v^\star \right\|^2\right] \\
			& \leq & p\alpha^2 \E\left[ \left\|\u^{n-1}\right\|^2\right] + (1-p)\E\left[\left\|\v^{n-1} - \v^\star \right\|^2\right] \\
			& = & \left(1-p\left(1-\alpha^2\right)\right)\E\left[\left\|\v^{n-1} - \v^\star\right\|^2\right] \label{eq-lem:dissipative-about-fixed-point-20}
		\end{eqnarray}
		where the first equality is due to \eqref{eq:asynchronous-identity}, the second equality is due to the translated system states, the third inequality is due to the $\alpha$-dissipative about $\v^\star$ assumption, and the final equality is due to substitution of the translated system states. Iterating the inequality \eqref{eq-lem:dissipative-about-fixed-point-20} $n$ times and taking a limit yields:
		\begin{eqnarray}
			\lim_{n \to \infty} \E\left[\left\|\v^n - \v^\star\right\|^2\right] & \leq & \left\|\v^0 - \v^\star\right\|^2 \lim_{n \to \infty} \left(1-p\left(1-\alpha^2\right)\right)^n \\
			& = & 0
		\end{eqnarray}
		where $1-p\left(1-\alpha^2\right)<1$ since $\alpha \in[0,1)$ by assumption and $p\in(0,1)$ by Definition~\ref{def:asynchronous}. We conclude this proof by showing that $\v^\star$ is the only fixed-point, i.e.~that $\mathcal{F}_{T}$ is singleton. Denote by $u^\star$ a different fixed-point of $T$ and let $\w = \u^\star - \v^\star$, then:
		\begin{equation}
			\|\w\| = \left\|\u^\star - \v^\star \right\| = \left\|T(\u^\star) - T(\v^\star) \right\| = \left\|T(\w+\v^\star) - T(\v^\star) \right\| \leq \alpha\left\| \w \right\| = \alpha \left\| \u^\star - \v^\star \right\| = \alpha \|\w\|
		\end{equation}
		which contradicts the assumption that $\v^\star \neq \u^\star$ since $\alpha < 1$ hence $\v^\star$ is the unique fixed-point of $T$.
	\end{pf}

	\begin{pf}[Lemma~\ref{lem:dissipative-about-some-point}]
		We begin in the synchronous implementation setting and show that for any $\epsilon > 0$ there exists a finite integer $n_0 \in \N$ such that the state evolution sequence $\left\{ \v^n \right\}_{n=n_0+1}^{\infty}$ is contained within the Euclidean ball centered at $\c$ with radius $\left\| T(\c) - \c\right\|(1-\alpha)^{-1}+ \epsilon$ irrespective of $\v^0$. Let $\un$ denote a translated by $\c$ state evolution sequence in $\Rk$, i.e. $\u^n = \v^n - \c$ for $n = 0,1,2,\dots$. We then obtain the inequality: 
		\begin{eqnarray}
			\left\| \v^n - \c \right\| & \leq & \left\| T(\v^{n-1}) - T(\c) \right\| + \left\| T(\c) - \c \right\| \\
									   & =    & \left\| T(\u^{n-1} + \c) - T(\c) \right\| + \left\| T(\c) - \c \right\| \\
									   & \leq & \alpha\left\| \v^{n-1} - \c \right\| + \left\| T(\c) - \c \right\| \label{eq-lem:dissipative-about-some-point-10}
		\end{eqnarray}
		where the first inequality is due to the triangle inequality, the second equality is due to the translated system states, and the third inequality is due to the assumed $\alpha$-connicity of $T$ at $\c$ and the translated system states. Iterating the inequality \eqref{eq-lem:dissipative-about-some-point-10} $n$ times yields:
		\begin{equation}
			\left\| \v^n - \c \right\| \leq \alpha^n\left\| \v^0 - \c \right\| + \left\| T(\c) - \c \right\| \sum_{k=0}^{n-1}\alpha^k  
		\end{equation}
		which we loosen to: 
		\begin{equation}
			\left\| \v^n - \c \right\| \leq \alpha^n\left\| \v^0 - \c \right\| + \frac{\left\| T(\c) - \c \right\|}{1-\alpha}.
		\end{equation}
		Therefore, taking $n_0 > \log_{\alpha}\epsilon - \log_{\alpha}\|\v^0 - \c\|$ is sufficient to ensure that $\alpha^n\left\| \v^0 - \c \right\|  < \epsilon$ for $n > n_0$. 

		In the asynchronous implementation setting we proceed to show that for every $\epsilon > 0$ there exists a finite integer $n_0\in\N$ such that the state evolution sequence $\left\{ \v^n \right\}_{n=n_0+1}^{\infty}$ also satisfies:
		\begin{eqnarray}
			\E\left[ \left\|\v^n - \c \right\| \right] \leq \frac{\left\| T(\c) - \c \right\|}{(1-\alpha)} + \epsilon.
		\end{eqnarray}
		To this end, we obtain the inequality:
		\begin{eqnarray}
			\E\left[ \left\| \v^n - \c \right\| \right] & = & p\E\left[ \left\| T(\v^{n-1}) - \c \right\| \right] + (1-p)\E\left[ \left\| \v^{n-1} - \c \right\|^2 \right] \\
			& \leq & p\E\left[ \left\| T(\v^{n-1}) - T(\c) \right\| \right] + p\left\| T(\c) - \c \right\| + (1-p)\E\left[ \left\| \v^{n-1} - \c \right\| \right] \\
			& \leq & \left(1 - p(1-\alpha)\right) \E\left[ \left\| \v^{n-1} - \c \right\| \right] + p\left\| T(\c) - \c \right\|  \label{eq-lem:dissipative-about-some-point-20}
		\end{eqnarray}
		where the first equality is due to \eqref{eq:asynchronous-identity} where we additionally replace $\v^\star$ by $\c$ since no fixed-point properties are utilized in obtaining \eqref{eq:asynchronous-identity}, the second inequality is due to the triangle inequality, and the third inequality is due to the translated state evolution sequence and the $\alpha$-connicity of $T$ at $\c$. Iterating the inequality \eqref{eq-lem:dissipative-about-some-point-20} $n$ times yields:
		\begin{equation}
			\E\left[ \left\| \v^n - \c \right\| \right] \leq \left(1-p(1-\alpha)\right)^n\left\| \v^0 - \c \right\| + p\left\| T(\c) - \c \right\|  \sum_{k=0}^{n-1}\left(1-p(1-\alpha)\right)^k
		\end{equation}
		which we loosen to: 
		\begin{equation}
			\E\left[ \left\| \v^n - \c \right\| \right] \leq \left(1-p(1-\alpha)\right)^n\left\| \v^0 - \c \right\| + \frac{p\left\| T(\c) - \c \right\|}{p(1-\alpha)}. 
		\end{equation}
		Therefore, there always exists a finite value of $n_0$ sufficiently large to ensure that $\left(1-p(1-\alpha)\right)^n\left\| \v^0 - \c \right\| < \epsilon$ for $n > n_0$ since $\lim_{n\to\infty}\left(1-p(1-\alpha)\right)^n\left\| \v^0 - \c \right\| = 0$. This concludes the proof.
	\end{pf}

	\begin{pf}[Lemma~\ref{lem:dissipative-everywhere}]
		We first derive two useful inequalities in the synchronous implementation setting. First, observe that substituting $\v' = \v^n$ and $\u' = \v^{n-1}$ into \eqref{eq:alpha-conic-everywhere-equivalent} and iterating the inequality $n-1$ times for an arbitrarily selected initial system state $\v^0$ gives:
		\begin{eqnarray} 
			\left\| \v^{n} - \v^{n-1} \right\| &\leq& \alpha \left\| \v^{n-1} - \v^{n-2} \right\| \\
			& \leq & \alpha^{n-1} \left\| \v^1 - \v^0 \right\|. \label{eq-lem:dissipative-everywhere-10}
		\end{eqnarray}
		Second, denote by $n$ and $m$ two positive integers satisfying $m < n$ and, again refering to \eqref{eq:alpha-conic-everywhere-equivalent} for an arbitrarily selected initial system state $\v^0$, assign $\v' = \v^n$ and $\u' = \v^m$. Then, by repeated application of the triangle inequality we obtain:
		\begin{eqnarray} \label{eq-lem:dissipative-everywhere-20}
			\left\| \v^n - \v^m \right\| \leq \sum_{k=m}^{n-1} \left\| \v^{k+1} - \v^{k} \right\|.
		\end{eqnarray}
		Together the inequlities in \eqref{eq-lem:dissipative-everywhere-10} and \eqref{eq-lem:dissipative-everywhere-20} imply that:
		\begin{eqnarray}
			\left\| \v^n - \v^m \right\| & \leq & \left( \sum_{k=m}^{n-1} \alpha^k \right) \left\| \v^1 - \v^0 \right\| \\
			& = & \frac{\alpha^m\left(1-\alpha^{n-m}\right)}{1 - \alpha} \left\| \v^1 - \v^0 \right\| \label{eq-lem:dissipative-everywhere-30}
		\end{eqnarray}
		from which the state evolution sequence $\vn$ is readily shown to be Cauchy. Subsequently, \eqref{eq-lem:dissipative-everywhere-30} establishes that the limit point $\v^\star$ exists and that the property $\v^n \to \v^\star$ is both well-defined and true. Moreover, from the continuity of the system operator $T$ implied by the $\alpha$-conic everywhere assumption it follows that:
		\begin{eqnarray}
			T(\v^\star) = T\left( \lim_{n \to \infty} \v^n \right) = \lim_{n \to \infty} T(\v^n) = \lim_{n\to \infty} v^{n+1} = \v^\star
		\end{eqnarray}
		thus $\v^\star$ is indeed a fixed-point of $T$. We now prove that the fixed-point $\v^\star$ is unique. Substituting $\v' = \v^\star$ and $\u' = \u^\star$ into \eqref{eq:alpha-conic-everywhere-equivalent} where $\v^\star$ and $\u^\star$ are assumed to be different fixed-points of $T$ yields:
		\begin{eqnarray}
			\left\| T\left(\v^\star\right) - T(\u^\star) \right\| \leq \alpha \left\| \v^\star - \u^\star \right\|
		\end{eqnarray}
		from which we conclude by contradiction that $\v^\star = \u^\star$ since $\alpha < 1$.

		We conclude by using the established existance and uniqueness of $\v^\star$ in proving that the state evolution sequence generated by an asynchronous implementation of $T$ converges to $\v^\star$. Indeed, such fixed-point properties belong to the system operator $T$ itself and are independent of the synchronous presentation details. Observe:
		\begin{eqnarray}
			\E\left[\|\v^n - \v^\star \|^2\right]& = & p\E\left[\left\|T(\v^{n-1}) - T(\v^\star) \right\|^2  \right] + (1-p)\E\left[\left\|\v^{n-1} - \v^\star \right\|^2  \right] \\ 
			& \leq & \left(1 + p(\alpha^2 - 1)\right)\E\left[\left\|\v^{n-1} - \v^\star\right\|^2 \right] \label{eq-lem:dissipative-everywhere-40}
		\end{eqnarray}
		where the equality is due to the \eqref{eq:alpha-conic-everywhere-equivalent} and the inequality is due to the $\alpha$-dissipativity everywhere property of $T$. Iterating the inequality \eqref{eq-lem:dissipative-everywhere-40} $n$ times and taking a limit gives:
		\begin{eqnarray}
			\lim_{n\to\infty} \E\left[\|\v^n - \v^\star \|^2\right] & \leq & \|\v^0 - \v^\star \|^2 \lim_{n \to \infty}(1+p(\alpha^2 - 1))^n \\ 
			& = & 0
		\end{eqnarray}
		where $\v^0$ is arbitrarily chosen and deterministic hence we have proven that $\v^n \to \v^\star$ in the sense of \eqref{eq:convergence-asynch}.
	\end{pf}

	\begin{pf}[Lemma~\ref{lem:dissipative-everywhere-ball}] We begin by proving a containment property of the system operator $T$ which will be critical to both the synchronous and asynchronous convergence results alike. For any element $\v \in B(\c,r)$ we have:
		\begin{eqnarray} \label{eq-lem:synch-dissipative-everywhere-ball-10}
			\left\| \c - T(\v) \right\| & \leq & \left\| \c - T(\c) \right\| + \left\| T(\c) - T(\v) \right\| \\
									    & \leq & (1-\alpha)r + \alpha\left\| \c - \v \right\|\\
									    & \leq & (1-\alpha)r + \alpha r  \\
									    &  =   & r  \label{eq-lem:synch-dissipative-everywhere-ball-20}
		\end{eqnarray}
		where the first inequality is due to the triangle inequality, the second inequality is due to the $\alpha$-dissipativity of $T$ over $B(\c,r)$, and the third inequality is due to the assumption that $\v\in B(\c,r)$. It immediately follows that $T$ cannot map any state $\v \in B(\c,r)$ outside of $B(\c,r)$. Therefore, for synchronous implementations we directly conclude that the state evolution sequence $\vn$ is forever contained within $B(\c,r)$ provided that $\v^0 \in B(\c,r)$.

		We now prove by induction that the state evolution sequence $\vn$ generated by an asynchronous implementation of $T$ is contained within $B(\c,r)$ in the sense of \eqref{eq:inside-ball} for an arbitrarily selected, deterministic initial system state  $\v^0$ inside $B(\c,r)$:
		\begin{itemize}
			\item[] Base case: For $n=0$ we have by assumption that $\v^0 \in B(\c,r)$.
			\item[] Induction step: Suppose that $\v^{n-1}$ is contained within $B(\c,r)$, i.e.~ $\E\left[\|\v^{n-1} -\c\|\right] \leq r$. We then obtain the inequality: 
			\begin{eqnarray} \label{eq-lem:synch-dissipative-everywhere-ball-30}
				\E\left[\left\| \v^n - \c \right\|\right] & \leq &\left(1 - p(1-\alpha)\right) \E\left[ \left\| \v^{n-1} - \c \right\| \right] + p\left\| T(\c) - \c \right\| \\
				& \leq & \left(1 - p(1-\alpha)\right) r + p(1-\alpha)r \\
				& = & r
			\end{eqnarray}
			where the first inequality is due to \eqref{eq-lem:dissipative-about-some-point-20}, the second inequality is due to the assumption of the induction step, and the third equality is due to the inequality \eqref{eq-lem:synch-dissipative-everywhere-ball-10}-\eqref{eq-lem:synch-dissipative-everywhere-ball-20}. Therefore, $\v^n$ is contained within $B(\c,r)$.
		\end{itemize}
		Said another way, on average a state evolution sequence which enters the ball $B(\c,r)$ never leaves. The remaining convergence claims for both synchronous and asynchronous implementation protocols follow by application of Lemma~\ref{lem:dissipative-everywhere} with the following assignments:  $B(\c,r)$, which is itself a bonafide complete metric space, takes the place of $\Rk$, the system operator $T$ is taken to be $T: B(\c,r) \to B(\c,r)$, and the initial system state $\v^0$ is arbitrarily selected from $B(\c,r)$.
	\end{pf}

	\begin{pf}[Lemma~\ref{lem:dissipative-everywhere-composition}]
		By application of Lemma~\ref{lem:dissipative-everywhere} we conclude that the system operator $T^m$ has a unique fixed-point $\u^\star$, i.e.~$T^m(\u^\star) = \u^\star$, since it is itself an $\alpha$-dissipative everywhere system operator. We proceed by arguing that $\u^\star$ is also a fixed-point of $T$ and thus $\u^\star = \v^\star$. Indeed, since $T^m$ is assumed to satisfy the functional translation property it follows that:
		\begin{equation}
			T^m\circ T(\u^{\star}) = T\circ T^m(\u^{\star})= T(\u^{\star}).
		\end{equation}
		Hence, by the uniqueness of the fixed-point of $T^m$ we have that $\u^\star = T(\u^\star)$ and therefore by definition $\v^\star = \u^\star$.

		We next show that the state evolution sequence $\vn$ associated with either a synchronous or asynchronous implementation protocol tends to $\v^\star$ and that the fixed-point $\v^\star$ is unique. In the synchronous setting, for each integer $i$ satisfying $0 \leq i \leq m - 1$, the $\alpha$-dissipative everywhere system operator $T^m$ generates a subsequence of $\vn$ consisting of the system states $\{\v^{mp+i}\}_{p=0}^{\infty}$, which is further illustrated by the relation:
		\begin{eqnarray}
			& \v^{mp+i} = \underbrace{T^m \circ T^m \circ \cdots \circ T^m }_{p\text{ compositions}}\circ T^i(\v^0), &  p \geq 0,
		\end{eqnarray}
		where $T^0$ is interpreted as the identity operator and $\v^0$ is an arbitrarily selected initial system state. Said another way, the state evolution sequence $\vn$ consists of the interleaved system states generated from $m$ synchronous implementations of $T^m$ with respective initial states $T^{i}(\v^0)$, $i=0,1,\dots,m-1$. Therefore, we conclude that: 
		\begin{eqnarray}
			\lim_{p \to \infty} \left\| \v^{mp+i} - \v^\star \right\| = 0, & & i = 0, 1, \dots, m-1.
		\end{eqnarray} 
		In the asynchronous setting, the proof follows the same argument. Namely, we note that the state evolution sequence $\vn$ generated using \eqref{eq:asynchronous} is decomposable again into $m$ subsequences each satisfying $\v^{mp+i} \to \v^\star$, i.e.~in the sense:
		\begin{eqnarray}
			\lim_{p \to \infty} \E\left[\left\| \v^{mp+i} - \v^\star \right\|^2 \right] = 0, & & i = 0, 1, \dots, m-1.
		\end{eqnarray} 
		Utilizing the fact that convergence for both synchronous and asynchronous implementation protocols is independent of $i$ and each subsequence converges to the same fixed-point $\v^\star$ we conclude that $\v^n \to \v^\star$ and that the system state $\v^\star$ is unique.
	\end{pf}

	\begin{pf}[Lemma~\ref{lem:dissipative-everywhere-filtered}]
		We prove this lemma for both synchronous and asynchronous implementation protocols by direct application of Lemma~\ref{lem:dissipative-everywhere}. To see this, observe that if the system operator $T$ is $\alpha$-dissipative everywhere then the filtered system operator $T_{f}(\rho, \cdot)$ is $\rho\alpha + |1-\rho|$-conic everywhere:
		\begin{eqnarray}
			\sup_{\u \neq \underline{0}} \frac{ \left\| T_{f}(\rho, \u+\v) - T_{f}(\rho, \v) \right\| }{ \|\u\| }
			& = & \sup_{\u \neq \underline{0}} \frac{ \left\| \rho \left[T(\u+\v) - T(\v)\right] + (1-\rho)\u \right\| }{ \|\u\| } \\
						& \leq & \sup_{\u \neq \underline{0}} \frac{   \left\|\rho\left[T(\u+\v) - T(\v)\right]\right\| }{ \|\u\| } + |1-\rho| \\
			& \leq & \rho\left(\sup_{\u \neq \underline{0}} \frac{   \left\|\left[T(\u+\v) - T(\v)\right]\right\| }{ \|\u\| }\right) + |1-\rho| \\ 
			& \leq & \rho\alpha + |1-\rho|
		\end{eqnarray}
		where the first equality is due to Definition~\ref{def:filtered-system}, the second inequality is due to the triangle inequality, the third inequality is due to the homogeneity of norms, and the fourth inequality is due to the $\alpha$-connic everywhere assumption on the system operator $T$. It is immediate that in order to ensure that $T_f$ is dissipative we require $\rho \in \left(0,\frac{2}{1+\alpha}\right)$ so that $\rho\alpha + |1-\rho|\in[0,1)$. The remaining claims follow by application of Lemma~\ref{lem:dissipative-everywhere} with the following assignments: $\rho T(\v) + (1-\rho)\v$ takes the place of $T(\v)$ and $\rho\alpha + |1-\rho|$ takes the place of $\alpha$. 
	\end{pf}

	\begin{pf}[Lemma~\ref{lem:passive-everywhere}]

		First we state two facts. 

		\begin{enumerate}
			\item[(i)] Let $\{a_k \colon k \in \N \}$ denote a bounded, non-negative sequence of real-valued scalars which additionally satisfy $\lim_{k\to\infty} ka_k > 0$, then $\sum_{k \in \N } a_k = \infty$.
			\item[(ii)] Let $T$ denote a passive everywhere system operator and define $f\colon \R \times \R \to \R$ as the scalar valued function 
		\begin{eqnarray} \label{lem:compact-set-T-10}
			f(l, u) \triangleq \inf_{\substack{l \leq \|\v - \v^\star \|^2 \leq u \\ \v^\star\in\mathcal{F}_{T}}} \left\| \v - T(\v)\right\|^2.
		\end{eqnarray}
		Then, $f(l,u) > 0$ for every pair of scalars $(l,u)$ satisfying $u>l> 0$.
 		\end{enumerate}
 		
 		Fact (i) is proven as follows. Let $a = \limsup_{k\to\infty} ka_{k} > 0$ and define $\{b_\ell \colon \ell\in\I\}$, $\I \subseteq \N$, as a subsequence of $\{ka_k \colon k \in \N\}$ such that $b_\ell \to a$. The existance of such a convergent subsequence is guaranteed, for example, by Theorem~3.17 in \cite{rudin}. The claim then follows directly from the inequality
		\begin{eqnarray}
			\sum_{k \in \N} a_k \geq \sum_{\ell \in \I}\frac{\ell a_\ell}{\ell} = \sum_{\ell \in \I} \frac{b_\ell}{\ell} = \infty
		\end{eqnarray}
		which is due to the subsequence definition and is tight (for partial sums) when $\I = \N$. 

		Fact (ii)  is a direct consequence of the fact that $T$ is continuous system operator and the objective function in \eqref{lem:compact-set-T-10} is defined over a non-empty compact set which does not contain any fixed-points of $T$.

		We now prove the theorem. Consider the state evolution sequence $\vn$ generated using the synchronous implementation protocol and observe that the sequence of scalars $\| \v^n - \v^\star\|^2$ for $n \geq 0$ is non-increasing: 
		\begin{eqnarray}
			\hspace{-0.3in} \left\| \v^{n} - \v^\star \right\|^2 \!\!\!\! & = & \!\!\! \rho \left\|  T(\v^{n-1}) - \v^\star \right\|^2 + (1-\rho)\left\| \v^{n-1} - \v^\star \right\|^2\! - \! \rho(1-\rho)\left\| T(\v^{n-1}) - \v^{n-1} \right\|^2 \\ 
			& \leq & \!\!\!\left\| \v^{n-1} - \v^\star \right\|^2 - \rho(1-\rho)\left\|T(\v^{n-1}) - \v^{n-1} \right\|^2 \label{eq-thm:sync-passive-inequality}
		\end{eqnarray}
		where the inequality is due to the passivity of $T$.  Iterating this inequality $n$ times yields
		\begin{eqnarray}
			\left\| \v^{n} - \v^\star \right\|^2 & \leq & \left\| \v^0 - \v^\star \right\|^2 - \rho(1-\rho)\sum_{m=0}^{n-1}\left\|T(\v^m) - \v^m \right\|^2 
		\end{eqnarray}
		and so the sequence $\left\| \v^{n} - \v^\star \right\|^2$ for $n\geq0$ is bounded above by $\|\v^0 - \v^\star\|^2$. Rearranging terms, utilizing the assumption that $\rho$ is restricted to the open unit interval, and loosening the inequality gives
		\begin{eqnarray}
			\sum_{m=0}^{n-1} \left\| T(\v^m) - \v^m \right\|^2 &\leq &\frac{1}{\rho(1-\rho)}\left\| \v^0 - \v^\star \right\|^2.
		\end{eqnarray}
		We conclude that the sequence $\|T(\v^n) - \v^n\|^2 \to 0$ by taking a limit (note the upper bound is independent of $n$) and applying the contrapositive to fact (i) above. The rate of this convergence also follows from this fact and is $o\left(\frac{1}{n}\right)$, i.e.~strictly faster than a $\frac{1}{n}$ sequence. In addition, since the state evolution sequence $\vn$ lies in the compact set $\{\v\in\Rk \colon \|\v - \v^\star\| \leq \|\v^0 - \v^\star\|\}$ it follows that the sequence has a limit point. We conclude that this limit point is an element $\v^\star \in \mathcal{F}_{T}$ since $T(\v^n) - \v^n \to \underline{0}$. This concludes the proof using the synchronous update mechanism.

		We now turn to the convergence of the state evolution sequence $\vn$ generated using an asynchronous implementation protocol, initially in an analagous manner to the presentation above. Indeed, observe that the sequence of scalars $\E\left[\left\| \v^n - \v^\star \right\|^2\right]$ for $n \geq 0$ is non-increasing:
		\begin{eqnarray}
			\E\left[\left\| \v^n - \v^\star \right\|^2\right] & = & p\E\left[\left\| T_f(\rho,\v^{n-1}) - \v^\star\right\|^2\right] +(1-p)\E\left[\left\|\v^{n-1} - \v^\star \right\|^2\right] \\
			& = & p\rho\E\left[\left\|T(\v^{n-1}) - \v^\star \right\|^2\right]  + p(1-\rho)\E\left[\left\|\v^{n-1} - \v^\star \right\|^2\right]  \\
			& & \enspace -p\rho(1-\rho)\E\left[\left\| T(\v^{n-1}) - \v^{n-1} \right\|^2\right] + p\E\left[\left\| \v^{n-1} - \v^\star \right\|^2\right] \nonumber   \\
			& \leq & \E\left[\left\|\v^{n-1} - \v^\star \right\|^2\right] - p\rho(1-\rho)\E\left[\left\| T\left(\v^{n-1}\right) - \v^{n-1}\right\|^2\right] 
		\end{eqnarray}
		where the first equality is due to \eqref{eq:asynchronous-identity} and the second equality is due \eqref{eq:pre-stewart} and the linearity of the expectation operator.  Iterating this inequality $n$ times results in the inequality
		\begin{eqnarray}
			\E\left[\left\| \v^{n} - \v^\star \right\|^2\right] & \leq & \left\| \v^0 - \v^\star \right\|^2 - p\rho(1-\rho)\sum_{m=0}^{n-1}\E\left[\left\|T(\v^m) - \v^m \right\|^2\right] \label{eq-app:upper-bound-async}
		\end{eqnarray}
		and so we conclude that the sequence $\E\left[\left\| \v^n - \v^\star \right\|^2\right]$ for $n \geq 0$ is bounded above by $\|\v^0 - \v^\star\|^2$. Rearranging terms, utilizing the assumption that $\rho$ is restricted to the open unit interval, and loosening the inequality results in
		\begin{eqnarray}
			\sum_{m=0}^{n-1} \E\left[\left\| T(\v^m) - \v^m \right\|^2\right] &\leq &\frac{1}{p\rho(1-\rho)}\left\| \v^0 - \v^\star \right\|^2.
		\end{eqnarray}
		Again, by application of the contrapositive of fact (i) above, we conclude that $\v^n - T(\v^n) \to \underline{0}$ in mean square (which further implies that $\v^n \to T(\v^n)$ in probability). We now show by contradiction that $\v^n \to \v^\star$ for some $\v^\star \in \mathcal{F}_{T}$ in mean square. Suppose that $ c = \lim_{n\to\infty}\E\left[\left\|\v^n - \v^\star \right\|^2\right]$ for some scalar $c > 0$. Note that the limit does indeed exist in part due to the  monotone and bounded properties of the sequence $\E\left[\left\| \v^n - \v^\star \right\|^2 \right]$. By application of the law of iterated expectation, we obtain the following inequality:
		\begin{eqnarray}\hspace{-.35in}
			\E\left[\left\| \v^n - \v^\star \right\|^2\right] & = & \E\left[\left\| \v^n - \v^\star \right\|^2
		  \mid \left\| \v^n - \v^\star \right\|^2 \geq \frac{c}{2}\right]\P\left(\left\| \v^n - \v^\star \right\|^2 \geq \frac{c}{2}\right) \\
			& & \enspace + \E\left[\left\| \v^n - \v^\star \right\|^2 \mid \left\| \v^n - \v^\star \right\|^2 < \frac{c}{2}\right]\P\left(\left\| \v^n - \v^\star \right\|^2 < \frac{c}{2}\right) \nonumber \\
			& \leq & \left\|\v^0 - \v^\star \right\|^2 \P\left(\left\| \v^n - \v^\star \right\|^2\geq \frac{c}{2}\right) + \frac{c}{2}\left(1-\P\left(\left\| \v^n - \v^\star \right\|^2 \geq \frac{c}{2}\right)\right)
		\end{eqnarray}
		where the upper bound for the term $\E\left[ \left\| \v^n - \v^\star \right\|^2 \mid  \left\| \v^n - \v^\star \right\|^2 \geq \frac{c}{2} \right]$ is from \eqref{eq-app:upper-bound-async}. Taking $n$ large enough and rearranging terms yields 
		\begin{eqnarray} \label{eq-thm:chain1}
 			0 < \frac{\frac{c}{2}}{\|\v^0 - \v^\star\|^2 - \frac{c}{2}} \leq \P\left(\left\| \v^n - \v^\star \right\|^2 \geq \frac{c}{2}\right).
		\end{eqnarray}
		By application of fact (ii) we are able to immediately conclude that 
		\begin{eqnarray} \label{eq-thm:chain2}
 			 \P\left(\left\| \v^n - \v^\star \right\|^2 \geq \frac{c}{2}\right) \leq \P\left( \left\| \v^n - T(\v^n) \right\|^2 \geq f\left(\frac{c}{2}, \|\v^0 - \v^\star\|^2\right) \right) 
		\end{eqnarray}
		since by the definition of the function $f$ in \eqref{lem:compact-set-T-10} the event $\|\v^n - \v^\star\|^2 \geq \frac{c}{2}$ implies $\|\v^n - T(\v^n)\|^2 \geq f(\frac{c}{2}, \|\v^0 - \v^\star\|^2)$. We then extend the chain of inequalities \eqref{eq-thm:chain1}-\eqref{eq-thm:chain2} by applying  Markov's inequality \cite{RVConvergence} to \eqref{eq-thm:chain2} and obtain
		\begin{eqnarray}
 			0 < \frac{\frac{c}{2}}{\|\v^0 - \v^\star\|^2 - \frac{c}{2}} \leq \P\left(\left\| \v^n - \v^\star \right\|^2 \geq \frac{c}{2}\right) \leq \frac{\E\left[\left\| \v^n -T(\v^n) \right\|^2\right]}{f\left(\frac{c}{2}, \|\v^0 - \v^\star\|^2 \right)}.
		\end{eqnarray}
		Taking a limit produces a contradiction since we have already proven that the term $\E\left[\left\| \v^n - \v^\star \right\|^2\right]$ goes to zero, therefore the scalar $c$ must be zero, i.e.  
		\begin{eqnarray}
			\lim_{n \to \infty} \E\left[ \left\| \v^n - \v^\star \right\|^2 \right] = c = 0.
		\end{eqnarray}
		This concludes the proof that the state evolution sequence converges to an element of $\mathcal{F}_{T}$ in mean square.
	\end{pf}

	\begin{pf}[Lemma~\ref{lem:expansive-everywhere}]
		We begin by proving that the state evolution sequence $\vn$ generated by a synchronous implementation of the filtered system operator $T_f$ converges to a limit point $\v^\star$ which is a fixed-point of $T$ using the restriction of the filtering coefficient indicated in \eqref{eq-lem:expansive-everywhere-filter-parameter} and the topological mixing property in \eqref{eq-lem:expansive-everywhere-mixing-condition}. To this end, we obtain the inequality: 
		\begin{eqnarray}
			\left\|\v^{n} - \v^\star \right\|^2 & = & \rho^2\left\| T(\v^{n-1}) - T(\v^\star) \right\|^2 + (1-\rho)^2\left\| \v^{n-1} - \v^\star \right\|^2 + 2\rho(1-\rho)\left\langle T(\v^{n-1}) - \v^\star, \v^{n-1} - \v^\star \right\rangle \\
			& \leq & \left( \alpha^2\rho^2 + (1-\rho)^2\right)\left\| \v^{n-1} - \v^\star \right\|^2 + 2\rho(1-\rho)\left\langle T(\v^{n-1}) - \v^\star, \v^{n-1} - \v^\star \right\rangle \\
			& \leq & \left( \alpha^2\rho^2 + (1-\rho)^2\right)\left\| \v^{n-1} - \v^\star \right\|^2 + 2\gamma\rho(1-\rho)\left\| T(\v^{n-1}) - T(\v^\star)\right\|\left\| \v^{n-1} - \v^\star \right\| \\
			& \leq &  \left( \alpha^2\rho^2 + (1-\rho)^2 + 2\gamma\alpha\rho(1-\rho)\right)\left\| \v^{n-1} - \v^\star \right\|^2\label{eq-lem:expansive-everywhere-10}
		\end{eqnarray}
		where the first equality is due to \eqref{eq-lem:passive-everywhere-20} and the definition of a fixed-point, the second inequality is due to the $\alpha$-expansivity of $T$, the third inequality is due to \eqref{eq-lem:expansive-everywhere-mixing-condition} and the definition of a fixed-point, and the final inequality is again due to the $\alpha$-expansivity of $T$. From \eqref{eq-lem:expansive-everywhere-10} we conclude that:
		\begin{eqnarray} \label{eq-lem:expansive-everywhere-20}
			\left\|\v^{n} - \v^\star \right\| & \leq &  \theta(\rho)^n \left\| \v^{0} - \v^\star \right\|.
		\end{eqnarray}
		where $\theta^2(\rho) =  \alpha^2\rho^2 + (1-\rho)^2 + 2\gamma\alpha\rho(1-\rho)$ is a convex quadratic form which we rewrite as:
		\begin{eqnarray} \label{eq-lem:expansive-everywhere-30}
			\theta^2(\rho) = \left(\alpha^2 + 1 - 2\alpha\gamma \right)\rho^2 + \left( 2\alpha\gamma - 2 \right)\rho + 1.
		\end{eqnarray} 
		We next show that $\theta^2(\rho)\in(0,1)$ for $\rho \in \left(0,\frac{2\left(1-\alpha\gamma\right)}{1+\alpha^2 - 2\alpha\gamma}\right)$. Note that $\theta^2(\rho)$ is guaranteed to be positive everywhere since it cannot have a real root. In particular, the condition:
		\begin{eqnarray}
			\left(2\alpha\gamma - 2  \right)^2 - 4\left(\alpha^2 + 1 - 2\alpha\gamma\right) < 0
		\end{eqnarray}
		is equivalent to the condition that $\gamma < 1$ which is ensured by the assumption in \eqref{eq-lem:expansive-everywhere-mixing-condition}. The minimum of $\theta^2(\rho)$ occurs at $\rho = \frac{\left(1-\alpha\gamma\right)}{1+\alpha^2 - 2\alpha\gamma}$ which is the midpoint of the interval in \eqref{eq-lem:expansive-everywhere-filter-parameter} and therefore the supremum of $\theta^2(\rho)$ is obtained at the limit points of this open interval. For any fixed-value of $\rho$ satisfying \eqref{eq-lem:expansive-everywhere-filter-parameter} verifying $\theta^2(\rho)$ involves straightforward algebra. We derive an upper bound for $\rho$ such that $\theta^2(\rho)$ is strictly upper bounded by unity by reducing the expression in \eqref{eq-lem:expansive-everywhere-30} by subtracting constants and dividing through by $\rho > 0$. This results in the condition $\rho(1+\alpha^2 - 2\alpha\gamma)  <  2(1-\alpha\gamma)$. By symmetry, the lower bound is given by $\rho > 0$. Therefore, we have shown that $\theta(\rho) \in (0,1)$ for
		\begin{eqnarray}  \label{eq-lem:expansive-everywhere-40}
			\rho \in \left(0, \frac{2(1-\alpha\gamma)}{1 + \alpha^2 - 2\alpha\gamma}\right).
		\end{eqnarray}
		Finally, we conclude that $\v^n \to \v^\star$ from \eqref{eq-lem:expansive-everywhere-20} by taking a limit:
		\begin{eqnarray}
			\lim_{n\to \infty} \left\| \v^{n} - \v^\star \right\| &\leq& \left\| \v^0 - \v^\star \right\| \lim_{n\to \infty}\theta(\rho)^{n} \\
			& = & 0.
		\end{eqnarray}
		The proof of convergence for state evolution sequence generated using an asynchronous implementation protocol follows analogously. Indeed, we begin by establishing the inequality:  
		\begin{eqnarray}`
			\E\left[ \left\| \v^n - \v^\star \right\|^2 \right] & = & p\E\left[ \left\| \rho(T(\v^{n-1})-\v^\star) + (1-\rho)(\v^{n-1}-\v^\star) \right\|^2 \right] + (1-p)\E\left[ \left\| \v^{n-1} - \v^\star \right\|^2 \right] \\
			& = & p\rho^2\E\left[ \left\| T(\v^{n-1})-T(\v^\star)\right\|^2 \right] + \left(1 - p + p(1-\rho)^2\right)\E\left[ \left\| \v^{n-1}-\v^\star \right\|^2 \right]  \\
			& & + 2p\rho(1-\rho)\E\left[\left\langle T(\v^{n-1}) - \v^\star, \v^{n-1} - \v^\star \right\rangle\right] \\
			& \leq &  \left(p\rho^2\alpha^2 + 1 - p + p(1-\rho)^2 + 2p\alpha\gamma\rho(1-\rho)\right)\E\left[ \left\| \v^{n-1})-\v^\star\right\|^2 \right] \label{eq-lem:expansive-everywhere-50}
		\end{eqnarray}
		where the second equality is due to the identity $\|\v+\u\|^2 = \langle \v + \u, \v + \u \rangle = \|\v\|^2 + \|\u\|^2 + 2\langle \v, \u \rangle$ and the definition of a fixed-point, and the third inequality is due to \eqref{eq-lem:expansive-everywhere-mixing-condition} and the $\alpha$-expansivity of $T$ and the mixing property in \eqref{eq-lem:expansive-everywhere-mixing-condition}. From \eqref{eq-lem:expansive-everywhere-50} we conclude that:
		\begin{eqnarray} \label{eq-lem:expansive-everywhere-70}
			\E\left[ \left\| \v^n - \v^\star \right\|^2 \right]  & \leq &  \phi(\rho)^n \left\| \v^{0} - \v^\star \right\|^2.
		\end{eqnarray}
		where $\v^0$ is an arbitrarily selected, determinisitc initial system state and $\phi(\rho) =  \alpha^2\rho^2 + (1-\rho)^2 + 2\gamma\alpha\rho(1-\rho)$ is a convex quadratic form which we rewrite as:
		\begin{eqnarray} \label{eq-lem:expansive-everywhere-80}
			\phi(\rho) = p\left(\alpha^2 + 1 - 2\alpha\gamma \right)\rho^2 + p\left( 2\alpha\gamma - 2 \right)\rho + 1.
		\end{eqnarray} 
		Note that $\phi(\rho) = \theta^2(\rho)$ for $p=1$. We next show that $\phi(\rho)\in(0,1)$ for $\rho$ satisfying \eqref{eq-lem:expansive-everywhere-filter-parameter}. $\phi(\rho)$ is guaranteed to be positive everywhere since it cannot have a real root. This condition in the asynchronous case is again equivalent to $\gamma < 1$. Indeed, the value of $\rho$ which minimizes $\phi(\rho)$ is easily seen to be the same as $\theta^2(\rho)$. We derive an upper bound for $\rho$ such that $\phi(\rho)$ is strictly upper bounded by unity by reducing the expression in \eqref{eq-lem:expansive-everywhere-80} by subtracting constants and dividing through by $p\rho > 0$. This again results in the condition $\rho(1+\alpha^2 - 2\alpha\gamma)  <  2(1-\alpha\gamma)$. By symmetry, the lower bound is again given by $\rho > 0$ hence we have for an asynchronous implementation protocol that $\phi(\rho) \in (0,1)$ provided that $\rho$ satisfies \eqref{eq-lem:expansive-everywhere-40}. Therefore, taking a limit of the inequality \eqref{eq-lem:expansive-everywhere-70} yields:
		\begin{eqnarray}
			\lim_{n\to \infty} \E\left[ \left\| \v^{n} - \v^\star \right\|^2 \right] &\leq& \left\| \v^0 - \v^\star \right\|^2 \lim_{n\to \infty}\phi(\rho)^{n} \\
			& = & 0
		\end{eqnarray}
		from which we conclude that $\v^n \to \v^\star$ in the sense of \eqref{eq:convergence-asynch}. 
	\end{pf}

\newpage
\FloatBarrier
\bibliographystyle{IEEEbib}
\nocite{*}
\bibliography{refs}

\end{document}